\documentclass[11pt]{article}

\usepackage[latin1]{inputenc}
\usepackage[T1]{fontenc}
\usepackage{amsmath,amssymb,amsthm}

\usepackage{latexsym}            

\DeclareMathOperator{\reg}{reg} \DeclareMathOperator{\depth}{depth}
\DeclareMathOperator{\im}{im} 
\DeclareMathOperator{\Supp}{Supp} \DeclareMathOperator{\HH}{H}
\DeclareMathOperator{\hho}{h}
\DeclareMathOperator{\Proj}{Proj} \DeclareMathOperator{\Spec}{Spec}
\DeclareMathOperator{\PGor}{PGor} 
\DeclareMathOperator{\CM}{CM} 
\DeclareMathOperator{\Hom}{Hom} \DeclareMathOperator{\Tor}{Tor}
\DeclareMathOperator{\Ext}{Ext} \DeclareMathOperator{\ext}{ext}
 
 \DeclareMathOperator{\Hilb}{Hilb}
\DeclareMathOperator{\GradAlg}{GradAlg}
 
\DeclareMathOperator{\Smc}{Smc} 
\DeclareMathOperator{\SmCM}{SmCM}
\DeclareMathOperator{\CI}{CI} \DeclareMathOperator{\CICM}{CICM}
 \DeclareMathOperator{\LevAlg}{LevAlg}
 \DeclareMathOperator{\length}{length}

\DeclareMathOperator{\Diff}{Diff} 

\newtheorem{theorem}{Theorem}
\newtheorem{lemma}[theorem]{Lemma}
\newtheorem{corollary}[theorem]{Corollary}
\newtheorem{definition}[theorem]{Definition}
\newtheorem{example}[theorem]{Example}
\newtheorem{proposition}[theorem]{Proposition}

\newtheorem{remark}[theorem]{Remark}

\newcommand{\pp}{{\mathbb P}}

\newcommand\sJ{{\mathcal J}}

\newcommand\sN{{\mathcal N}}
\newcommand\sO{{\mathcal O}}

\newcommand{\proj}[1]
{ \mathchoice
           { {\mathbb P}^{#1} }
           { {\mathbb P}^{#1} }
           { {\mathbb P}^{#1} }
           { {\mathbb P}^{#1} }
         }

\begin{document}
\bibliographystyle{plain}
\title{Families of Artinian and one-dimensional algebras}
\author{Jan O. Kleppe\thanks{Partially
  supported by Research Council of Norway, project no. 154077/420.}}
\date{\ }
\maketitle
\vspace*{-0.75in}
\begin{abstract}
\noindent The purpose of this paper is to study families of Artinian or one
dimensional quotients of a polynomial ring $R$ with a special look to level
algebras. Let $\GradAlg^H(R)$ be the scheme parametrizing graded quotients of
$R$ with Hilbert function $H$. Let $B \rightarrow A$ be any graded surjection
of quotients of $R$ with Hilbert function $H_B$ and $H_A$, and $h$-vectors
$h_B=(1,h_1,...,h_j,...)$ and $h_A$, respectively. If $\depth A = \dim A \leq
1$ and $A$ is a ``truncation'' of $B$ in the sense that
$h_A=(1,h_1,...,h_{j-1},\alpha,0,0,...)$ for some $\alpha \leq h_j$, then we
show there is a close relationship between $\GradAlg^{H_A}(R)$ and
$\GradAlg^{H_B}(R)$ concerning e.g.\! smoothness and dimension at the points
$(A)$ and $(B)$ respectively, provided $B$ is a complete intersection {\rm or}
provided the Castelnuovo-Mumford regularity of $A$ is at least 3 (sometimes 2)
larger than the regularity of $B$. In the complete intersection case we
generalize this relationship to ``non-truncated'' Artinian algebras $A$ which
are compressed or close to being compressed. For more general Artinian
algebras we describe the dual of the tangent and obstruction space of
deformations in a manageable form which we make rather explicit for level
algebras of Cohen-Macaulay type 2. This description and a linkage theorem for
families allow us to prove a conjecture of Iarrobino on the existence of at
least two irreducible components of $\GradAlg^H(R)$, $H=(1,3,6,10,14,10,6,2)$,
whose general elements are Artinian level algebras of type 2.

\noindent {\bf AMS Subject Classification.} 14C05, 13D10, 13D03,
13D07, 13C40, 13D02.

\noindent {\bf Keywords}.  Parametrization, Artinian
algebra, level algebra, Gorenstein algebra, licci, Hilbert scheme,
duality, algebra (co)homology, canonical module, normal module.
 \end{abstract}
 \vspace*{-0.25in}
\thispagestyle{empty}

\section{Introduction}
The main goal of this paper is to contribute to the classification of
Artinian and one dimensional graded quotients of a polynomial ring $R$
in $n$ variables (of degree one) over an algebraically closed field
$k$. In particular we study the scheme $ \GradAlg^{H}(R)= \GradAlg(H)
$ which parametrizes graded quotients $A$ of $R$ of $\depth A \geq
\min(1,\dim A)$ and with Hilbert function $H$.  $\GradAlg^{H}(R)$ is
the representing object of a correspondingly defined functor of flat
families and it may be non-reduced. Thus $ \GradAlg^{H}(R)$ may be
different from the parameter spaces studied by Iarrobino, Gotzmann and
others which study the ``same'' scheme with the reduced scheme
structure. In our approach we try to benefit of having a well
described tangent and obstruction space of $ \GradAlg^{H}(R)$ at $(A)$
at our disposal.

An important technique in determining $ \GradAlg(H_A)$ is to take a
graded surjection $B \rightarrow A$ of quotients of $R$ with Hilbert
functions $H_B$ and $H_A$ respectively, and, under certain conditions,
make the relationship between $\GradAlg(H_A)$ and $ \GradAlg(H_B)$
as tight as possible. We review some results of this technique in
Section 1. If $B = R/I_B$, let $N_B:= \Hom_R(I_B,B)$ and let
$\reg(B)=\reg(I_B)-1$ be the Castelnuovo-Mumford regularity of $B$.
Let
\begin{equation*} 
  ... \rightarrow  \oplus_{i=1}^{r_2} R(-n_{2,i}) 
  \rightarrow  \oplus_{i=1}^{r_1} R(-n_{1,i})  \rightarrow R \rightarrow B
  \rightarrow 0  \hspace{1 cm} 
\end{equation*} 
be the minimal resolution and let $\epsilon(A/B) = \sum_{i=1}^{r_1}
[H_B({n_{1,i}}) -H_A({n_{1,i}})]$. Our main results in Section 2
apply to $ \GradAlg(H_A)$ where $A$ is one-dimensional. We prove

\begin{theorem} \label{intromainzero}
  Let $R$ be a polynomial $k$-algebra and let $B = R/I_B
  \rightarrow A = R/I_A$ be a graded morphism such that $A$ is
  Cohen-Macaulay of dimension one and $\depth_{\mathfrak{m}}B \geq 1$
  and suppose $X:= \Proj(A) \hookrightarrow Y:=\Proj(B)$ is a local
  complete intersection of codimension $r \geq 0$. Let $H_A(v)= s$ for
  $v >>0$ and suppose either
  
  (a) \ \ \ \ $I_B$ is generated by a regular sequence (allowing $R = B$), or
  
  (b) \ \ \ \ $B_v \rightarrow A_v$ is an isomorphism for all $v \leq
  \max_{i}\{n_{2,i}\}$ and $\dim R - \dim B \geq 2$. \\
  Moreover suppose there is an integer $j$ such that $B_v \simeq A_v$ for all
  $v \leq j-1$ and such that $I_A$ is $(j+1)$-regular (i.e.  $\reg(A) \leq
  j$). Then $ \ \dim (N_A)_0= \dim (N_B)_0 + rs - \epsilon(A/B) \ , $ and $$
  \ 
  \dim_{(A)} \GradAlg^{H_A}(R)= \dim_{(B)} \GradAlg^{H_B}(R)+rs- \epsilon(A/B)
  \ .$$
  In particular $A$ is unobstructed as a graded $R$-algebra (i.e. $
  \GradAlg^{H_A}(R)$ is smooth at $(A)$) if and only if $B$ is unobstructed as
  a graded $R$-algebra.
\end{theorem}

One may look upon the conditions on $j$ above as assuming the minimal free
resolution of $I_{A/B}:=I_A/I_B$ to be semi-linear (close to being
linear, cf. \eqref{linear}), and the condition of (b) and (a) as
requiring this $j$ to be large enough in the case $Y$ is not a
complete intersection (CI).

Theorem~\ref{intromainzero} is  what we need to treat the
case where $X$ consists of $s$ points in generic position on $Y$ (so $H_A$ is
the truncation of $H_B$ at the level $s$). Indeed Geramita et al.
(\cite{GMR}) defines such a truncated Hilbert function by
$$H_A(i) = \inf \{H_B(i), s\} \ , \ $$
and they show that there exists a
reduced scheme $X$ on $Y$ with truncated Hilbert function $H_A$ provided $Y$
is reduced and consists of more that $s$ points. We prove

 \begin{corollary} \label{introsgeneric} Let $Y=\Proj(B)$, $B=R/I_B$, be a
   reduced scheme consisting of more than $s$ points, and let
   $X=\Proj(A)$ be $s$ points (avoid $SingY$) of codimension $r$ in
   generic position on $Y$.  Let $j$ be the smallest number such that
   $H_A(j) \neq H_B(j)$. If $Y$ is not a CI, suppose $j \geq
   \reg(I_B)+2$. Then $ \ \dim (N_A)_0= \dim (N_B)_0 + rs -
   \epsilon(A/B) \ , $ and $$
   \ \dim_{(A)} \GradAlg^{H_A}(R)=
   \dim_{(B)} \GradAlg^{H_B}(R)+rs- \epsilon(A/B) \ .$$
   Hence $A$ is
   unobstructed as a graded $R$-algebra iff $B$ is unobstructed as a
   graded $R$-algebra.
\end{corollary}

Moreover in Corollary~\ref{introsgeneric} (and Theorem~\ref{intromainzero}) we
may allow the codimension $r$ to vary along the $s$ points, say such that the
$i$-th point has codimension $c_i$ in $Y$ ($Y$ need not be equidimensional).
Then Corollary~\ref{introsgeneric} holds if we replace $rs$ by $\sum_{i=1}^s
c_i$.

The analogue of Theorem~\ref{intromainzero} for Artinian algebras is
the main result of Section 3.

\begin{theorem} \label{intromainartin} Let $R$ be a polynomial
  $k$-algebra and let $B = R/I_B \rightarrow A = R/I_A$ be a
  graded morphism such that $A$ is Artinian and $\ 
  \depth_{\mathfrak{m}}B \geq \min(1,\dim B)$, and suppose either
  
  (a) \ \ \ \ $I_B$ is generated by a regular sequence (allowing $R = B$), or
  
  (b) \ \ \ \ $B_v \rightarrow A_v$ is an isomorphism for all $v \leq
  \max_{i}\{n_{2,i}\}$ and $\dim R - \dim B \geq 2$. \\ 
  Let $F$ be a free $B$-module such that $F \rightarrow I_{A/B}$ is surjective
  and minimal, and suppose there is an integer $j$ such that the degrees of
  minimal generators of the $B$-module $\ker(F \rightarrow I_{A/B})$ $> j$
  (e.g. $B_v \simeq A_v$ for all $v \leq j-1$) and such that $I_A$ is
  $(j+1)$-regular (i.e.  $A_{j+1} = 0$).  Then $ \ \dim (N_A)_0= \dim (N_B)_0
  + \ \dim\, _0\! \Hom_B(F,A) - \epsilon(A/B) \ , $ and
 $$
  \ \dim_{(A)} \GradAlg^{H_A}(R)= \dim_{(B)} \GradAlg^{H_B}(R)+ \
  \dim\, _0\! \Hom_B(F,A) - \epsilon(A/B) \ .$$
  In particular $A$ is
  unobstructed as a graded $R$-algebra if and only if $B$ is
  unobstructed as a graded $R$-algebra.
\end{theorem}

In Proposition~\ref{corart} we show an improvement of
Theorem~\ref{intromainartin}(a) in the case ``$B_v \simeq A_v$ for
all $v \leq j-1$''. In this case we can skip the condition
$A_{j+1}=0$, or equivalently $(K_A)_{-j-1}=0$ ($K_A$ canonical module)
provided the minimal resolution of $K_A$ had no relations in degree
greater or equal to $j$. This generalization applies to algebras which
are compressed or close to being compressed. In the compressed case
the dimension and the smoothness of $\GradAlg(H_A)$ coincide with the
results of \cite{I84}.

Theorem~\ref{intromainartin} applies nicely to Artinian truncations
and more generally to Artinian quotients $A$ with $h$-vector
$H_A=(1,h_1, h_2,...,h_{j-1}, \alpha,0,0,..)$ where
$H_B=(1,h_1,h_2,...,h_{j-1},h_{j},h_{j+1},...)$ and $\alpha \leq h_j$.
In that case the relationship between $ \GradAlg(H_A)$ and the open
subset $ \GradAlg(H_B)_n$ of $ \GradAlg(H_B)$ consisting of points
$(B)$ where 
$\reg(I_B) \leq n$, may be described by an incidence correspondence
\begin{equation} \label{introincidartin}
  \begin{array}[h]{ccc}
    \GradAlg(H_B,H_A)_n & \stackrel{q}{\longrightarrow} & \GradAlg(H_B)_n
    \subset \GradAlg(H_B) \\
    \rule[-2mm]{0pt}{6mm} \downarrow ^{p} \\
    \GradAlg(H_A)
  \end{array}
\end{equation}
where $p$ and $q$ are the natural projections (cf. \eqref{incid} for details).

\begin{proposition} \label{introsgenartin} Let $H_B=(1,h_1,h_2, ...)$ be
  the Hilbert function of an algebra $B$ satisfying $\depth_{\mathfrak{m}}B
  \geq 1$, and let $j$, $n \leq j-2$ and $\alpha \leq h_j$ be integers. Let
  $H_A=(1,h_1,...,h_{j-1}, \alpha,0,0,..)$ and look to the maps $p$ and $q$ in
  \eqref{introincidartin}.  Then
  
  (i) \ \ $q$ is smooth and surjective with connected fibers, of fiber
  dimension $\alpha(h_j - \alpha)$, and
  
  (ii) \ \ $p$ is an isomorphism onto an open subscheme of $\ 
  \GradAlg(H_A)$. \\
  In particular the incidence correspondence \eqref{introincidartin}
  determines a well-defined injective application $\pi$ from the set of
  irreducible components $W$ of $ \GradAlg(H_B)_n$, to the set of irreducible
  components $V$ of $\ \GradAlg(H_A)$ whose general elements satify the Weak
  Lefschetz property. In this application the generically smooth components
  correspond.  Indeed $V=\pi(W)$ is the closure of $p(q^{-1}(W))$, and we have
  $$
  \ \dim V = \ \dim W + \  \alpha(h_j - \alpha) \ .$$
\end{proposition}

Also Theorem~\ref{intromainzero} allows a corollary very similar to
Proposition~\ref{introsgenartin} in the one dimensional case (cf.
Proposition~\ref{sgen}). 

In Section 4 we characterize the tangent and obstruction space of $
\GradAlg^{H}(R)$ at an Artinian algebra $(A)$. Note that if $A$ is
Gorenstein with socle degree $j$, then the obstruction space is the
dual of the kernel of the natural map $(S_2I_A)_j \rightarrow
({I_A}^2)_j$, or equivalently, the cokernel of $(\Lambda^2 I_A)_j
\rightarrow (I_A \otimes I_A)_j \simeq \Tor^R_1(I_A,K_A)_0$. This
result generalizes to the following result, in which $ \HH_2(R,A,K_A)$
is the algebra homology, cf. \eqref{alghom}.

\begin{theorem}\label{intromainartingrad}
  Let $R \rightarrow A = R/I_A$ be a graded Artinian quotient with
  Hilbert function $H$. Then $\dim (I_A \otimes_R K_A)_0$ is the
  dimension of the tangent space of $ \GradAlg^{H}(R)$ at $(A)$, and
  the dual of $ _0\!\HH_2(R,A,K_A)$ contains the obstructions of
  deforming $A$ as a graded $R$-algebra. In particular $
  \GradAlg^{H}(R)$ is smooth at $(A)$ 
  provided the natural ``antisymmetrization'' map $$I_A \otimes_R I_A
  \otimes_R K_A \rightarrow \Tor^R_1(I_A,K_A)$$
 (cf. \eqref{TorK}) is surjective in
  degree zero. 
\end{theorem}
Since we in Theorem~\ref{mainlev} show that the parameter space
$\L(H)$ of level algebras, introduced in \cite{GC} through the ideas
of \cite{IK}, essentially is an open subscheme of $ \GradAlg^{H}(R)$,
we get that Theorem~\ref{intromainartingrad} holds if we everywhere
replace $\GradAlg^{H}(R)$ by $\L(H)$ at a level algebra $(A)$. Note
that the tangent space of $\L(H)$ is already well described in
\cite{GC}.

Finally we look to type 2 level algebras $A=R/ann(F_1,F_2)$ where
$F_1$ and $F_2$ are forms of the same degree $s$ in the ``dual''
polynomial algebra of $R$. Such algebras are studied in \cite{I04},
and an extended draft of \cite{I04} determines the tangent space of $
\GradAlg^{H}(R)$ at $(A)$. If $H_A(i)= \min \{\dim R_i,
H_{A_1}(i)+H_{A_2}(i) \}$ for any $i$, then $\{F_1,F_2\}$ is said to
be complementary \cite{I04}. Using Theorem~\ref{intromainartingrad} we
describe the tangent and obstruction space of $ \GradAlg^{H}(R)$ at
$(A)$ in the following way.
 
\begin{proposition}\label{intromainlev2}
  Let $\{F_1,F_2\}$ be complementary forms of degree $s$, and let $A =
  R/I_A$ be the Artinian level quotient with Hilbert function $H$
  given by $I_A=ann(F_1,F_2)$. Let $I_{A_i}=ann(F_i)$.  Then $(I_A/I_A
  \cdot I_{A_1})_s \oplus (I_A/I_A \cdot I_{A_1})_s$ is the dual of
  the tangent space of $ \GradAlg^{H}(R)$ at $(A)$, and $
  _s\!\HH_2(R,A,A_1) \oplus \ _s\!\HH_2(R,A,A_2) $ is the dual of a
  space containing the obstructions of deforming $A$ as a graded
  $R$-algebra.  In particular if the sequences
\begin{equation*}  
   I_A \otimes_R I_A  \stackrel{\lambda}{\longrightarrow} I_A
   \otimes I_{A_i} \rightarrow I_A \cdot I_{A_i}   
   \end{equation*}
     where $ \lambda(x \otimes y)= x \otimes y - y \otimes x$, are exact for
   $i=1$ and $2$, then $\GradAlg(H_A)$ is unobstructed at $(A)$ and we have
  $
  \ \dim_{(A)} \GradAlg^{H_A}(R)= \ \sum_{i=1}^2 \dim (I_A/I_A \cdot
  I_{A_i})_s \ .$ 
\end{proposition}

Then we use Proposition~\ref{intromainlev2} and a linkage theorem
(Theorem~\ref{mainlink}) to prove a conjecture of A. Iarrobino,
appearing in the draft of \cite{I04}, namely that $\L(H)$ with
$H=(1,3,6,10,14,10,6,2)$ contains at least two irreducible components
whose general elements are level quotients of type 2
(Example~\ref{lev2comp}). Once having one example of such a phenomena,
we produce infinitely many by liaison (Remark~\ref{remlev2comp}). Even
though this conjecture was open until now, Iarrobino and Boij have in
a joint work already constructed other examples of reducible $\L(H)$
whose general elements are type 2 level quotients, one with
$H=(1,3,6,10,14,18,20,20,12,6,2)$, and seems to have got a doubly
infinite series of such components.

In this paper we give many examples to illustrate our results, some of
them with the use of Macaulay 2. Among examples of particular
interest, in addition to the proven conjecture, we mention
Example~\ref{twocomp} of two irreducible components of $\GradAlg(H)$
(and of $\PGor(H)$) whose intersection contains Artinian Gorenstein
algebras, and Example~\ref{exlink} of two components of
$\GradAlg(H)$ with $H=(1,3,6,6,3,1)$ whose general elements are
Artinian licci algebras.

I thank Anthony Iarrobino for interesting discussion on the subject
and for introducing me to type 2 level algebras and his conjecture by
giving me the extended draft of \cite{I04}. I also thank him and the
Northeastern University for their hospitality during my visit to
Boston in February 2005.
%
 
\subsection{Preliminaries}

Let $R$ be a polynomial $k$-algebra in $n$ variables of degree $1$
where $k$ is algebraically closed. In the following we focus on the
scheme parametrizing Artinian graded quotients $B$ of $R$, as well as
closed schemes in $\pp = \pp^{n-1}$, with fixed Hilbert function $H$.
Both schemes are denoted by $\GradAlg^{H}(R)$. If $H(v):=\dim B_v$
does not vanish for $v >> 0$, we call it the postulation Hilbert
scheme because this name seems to be most common, at least when it is
endowed with its reduced scheme structure and $\dim B =1$ (cf.\!
\cite{Go}, \cite{IK}). Since $\GradAlg^{H}(R)$ is the representing
object of a certain
functor of flat deformations, it may be non-reduced. 
We continue denoting it by $\GradAlg^{H}(R)$, to make it clear that it
may be non-reduced.
 
Now we recall the definition of $ \GradAlg^H(R)$. Let $ \Hilb^{ p} (\pp)$ be
the Hilbert scheme (\cite{SB}) parametrizing closed subschemes
$Y$ of $\pp=\Proj R$ with Hilbert polynomial $p \in \mathbb Q [t]$. 
The $k$-point of $\!\ \Hilb^{ p} (\pp)$ which corresponds to the $Y$ is
denoted by $(Y)$. A closed subscheme $Y$ of $\pp$ is called {\em
  unobstructed} if $\!\ \Hilb^{p}(\pp)$ is smooth at $(Y)$.

Let $\GradAlg(H):= \GradAlg^H(R)$ be the stratum of $\Hilb^{p}(\pp)$
given by deforming $Y \subset \pp$ with constant Hilbert function $
H_Y = H$ (i.e. its functor deforms the {\em homogeneous coordinate
  ring}, $B=R/I_B$, of $Y$ flatly), cf. \cite{K98} or
\cite{K04}. $\GradAlg^{H}(R)$ allows a natural scheme structure whose
tangent (resp.\ ``obstruction'') space at $(Y)$ is $\ 
_0\!\Hom_B(I_{B}/I_B^2,B) \simeq \ _0\!\Hom_R(I_{B},B) \ $ (resp.\ 
$_0\!\HH^2 (R,B,B)$), i.e.\ it is given by deforming $B$ as a graded
$R$-algebra (\cite{K79}, Thm.\! 1.5). In the case $H(v)$ does not
vanish for large $v$ (i.e. $B$ is non-Artinian), we may look upon
$\GradAlg^{H}(R)$ as parametrizing graded $R$-quotients, $R
\rightarrow B$, satisfying $\depth_{\mathfrak m}B \geq 1$ and with
Hilbert function $H_B = H$. If $B$ is Artinian, $\GradAlg^{H}(R)$
still represents a functor parametrizing graded $R$-quotients with
Hilbert function $H_B = H$ (see \cite{K04}, Prop.\! 9 and Thm.\! 11).
$B$ is called unobstructed as a graded $R$-algebra if and only if
(iff) $\GradAlg(H_B)$ is smooth at $(B)$, i.e. at $(Y)$.

\begin{remark} \label{conn}
  
  (a) It follows from a theorem of Pardue (Thm.\! 34 of \cite{Par}, cf.
  \cite{Go} for the codimension $2$ case) that $\GradAlg(H)$ is a connected
  scheme (see also \cite{Ma}).
  
  (b) One may use the representability of the functor which defines
  $\GradAlg(H)$ (through flat deformations) and the semicontinuity of
  the graded Betti numbers (i.e. of the number of minimal generators
  of a finitely generated module over $R$) to generalize
  Ragusa-Zappala's result for zero-schemes (\cite{RZ}), that different
  minima of the set of graded Betti numbers corresponds to different
  components of $\GradAlg(H)$, to any $H$ and $R$. Thus incomparable
  sets of graded Betti numbers lead to different components in general
  (see \cite{Mi} for a discussion).
\end{remark}

 The following comparison result is due to Ellingsrud
  (\cite{El}) in the case $\depth_{\mathfrak m}B \geq 2$, see \cite{K79},
  Thm.\! 3.6 and Rem.\! 3.7 for the general case. Below $s(I_B)$ is the
  minimal degree of the minimal generators of $I_B$. Note that the openness
  statements follow easily from the first isomorphism by the semicontinuity
  of $\ \dim \HH^1(Y,\widetilde {I_B}(v))$.

\begin{proposition} \label{Grad} Let $B = R/I_B$ satisfy
  $\depth_{\mathfrak m}B \geq 1$ and let $Y = \Proj(B)$. Then
  $$\GradAlg^{H}(R) \simeq \Hilb^{ p} (\pp) \ \ \ {\rm at} \ \ \ (Y) \ 
  ,$$
  provided $\ _0\!\Hom_R (I_B,\HH_{\mathfrak m}^1(B)) = 0$ (e.g.
  provided $\depth_{\mathfrak m}B \geq 2$). In particular the open
  sets $$U(H):= \{(B) \in \GradAlg^{H}(R) \arrowvert \ 
  _v\!\HH_{\mathfrak m}^1(B) = 0 \ {\rm for\ every \ } v \geq s(I_B)
  \}
  $$
  and $\ \{(B) \in \GradAlg^{H}(R) \arrowvert \depth_{\mathfrak m}B \geq 2
  \} $ of $\ \GradAlg^{H}(R)$ are also open in $\Hilb^{p}(\pp)$.
\end{proposition}
 Here
$\depth_{\mathfrak m}{M}$ ($M$ finitely generated) denotes the length of a
maximal $M$-sequence in the irrelevant maximal ideal $\mathfrak m$, and
$\HH^i_{\mathfrak m}(-)$ is the right derived functor of the functor
of sections with support in $\Spec (B/{\mathfrak m})$. Note that
$\depth_{\mathfrak m} M \geq r$ iff $\HH^i_{\mathfrak m}({M})=0$ for $i<r$, cf.
\cite{HARLOC}.
%
A Cohen-Macaulay $B$-module $M$ satisfies $\depth {M}= \dim {M}$ by
definition. If $B$ is Cohen-Macaulay of codimension $c$ in $R$ and
$K_{B}= \Ext^{c}_{R}(B,R(-n))$ is the canonical module of $B$, we know
by Gorenstein duality that $v$-graded piece of $\HH_{\mathfrak m}^i({
  M})$ satisfies
$$_{v}\!\HH_{\mathfrak{m}}^{i}(M) \simeq \ _{-v}\!\Ext_B^{n-c-i}(M,K_B)^{\vee} 
\ .$$

Two graded quotients, ${R}/{J}$ and ${R}/{J'}$, are said to be
linked by a complete intersection if there exists a homogeneous complete
intersection ideal ${L}$ such that ${J}={L}:{J'}$ and ${J'}={ L}:{J}$ (with
${L} \subseteq {J} \cap {J'})$.  The relationship of being linked generates
the equivalence relation, ``linkage''. ${B}={R}/{I_B}$ is said to be licci
(and hence Cohen-Macaulay) if it is in the linkage class of a complete
intersection (cf. \cite{MIG} for a survey). 
 
The algebraic (co)homology groups $\HH_2(R,B,M)$ and $\HH^2(R,B,M)$ may be
described as follows. The former group is given by an exact sequence
\begin{equation} \label{alghom}
   0 \rightarrow \HH_2(R,B,M) \rightarrow  {\HH}_1 \otimes_B M \rightarrow
  {G}_{1} \otimes_{R} {M} \rightarrow {I_B}/{I^2_B} \otimes_B M \rightarrow 0.
\end{equation}
in which ${G}_1$ is $R$-free, ${G}_1 \twoheadrightarrow {I_B}$ is surjective
and minimal, and $\HH_1=\HH_1(I_B)$ is the degree one Koszul homology of $I_B$
\cite{LS}. 
For the graded group $\HH^2(R,B,M)$ we only remark that by \cite{AND},
Prop.\! 16.1, and \cite{LS}, there are injections
\begin{equation} \label{cohom}
     _0\!\Ext_B^1(I_B/I_B^2,M)  \hookrightarrow \ _0\!\HH^2(R,B,M)
     \hookrightarrow \  _0\!\Ext_R^1(I_B,M)  
 \end{equation}  

 A quotient $B=R/I_B$ of codimension $c:= \dim R - \dim B$ in $R$ has a
 minimal $R$-free resolution of the following form (cf. \cite{EIS})
\begin{equation} \label{resB}
  ... \rightarrow {G}_c \rightarrow . . . 
  \rightarrow {G}_1 \rightarrow R \rightarrow B  \rightarrow 0  \hspace{0.2cm}
  ,  \hspace{1 cm} G_j = \oplus_{i=1}^{r_j} R(-n_{j,i}) 
\end{equation}
and $B$ is Cohen-Macaulay (CM) iff $G_{c+1}=0$. The function
$\max_{i}\{n_{j,i}\}-j$ is increasing as a function in $j$ if $B$ is
CM. If $B$ is Artinian (i.e. $c=n$), then $\max_{i}\{n_{c,i}\}-c$ is
the socle degree of $B$. More generally the Castelnuovo-Mumford
regularity of $I_B$ is given by $ \reg(I_B)=
\max_{\{j,i\}}{\{n_{j,i}-j+1}\}$ and $\reg(B) =\reg(I_B)-1$ (cf.
\cite{MIG}, p.\! 8). In particular
\begin{equation*} \label{Cast}
  \max_{i}\{n_{j,i}\}  \leq  \reg(I_B) +j-1 \ \ {\rm \ for \ any} \ j.
\end{equation*}

If $G_{c+1}=0$ and $G_c$ has rank $1$ (resp. $G_c = R(-s)^t$), then
$B$ is Gorenstein (resp. level of type $t$). In these cases $B$ is a
compressed Artinian $R$-algebra if $H_B$ (i.e. $H_B(v)$ for any $v$)
is as large as possible for a fixed socle degree and fixed type (cf.
\cite{I84} for existence).

An $R$-module $M$ of projective dimension $t-1$ is said to have a {\bf
  semi-linear} (resp. linear) resolution provided the minimal
resolution of $M$ has the following form
\begin{equation} \label{linear}
  0 \rightarrow R(-j-t)^{\beta_t} \oplus  R(-j-t+1)^{\alpha_t} \rightarrow . .
  .  \rightarrow 
  R(-j-1)^{\beta_1} \oplus  R(-j)^{\alpha_1} \rightarrow M  \rightarrow 0    
\end{equation}
(resp. with $\alpha_i=0$ for any $i$). With $B$ as in \eqref{resB} and $B
\rightarrow A \simeq B/I_{A/B}$ a graded surjection, we define
\begin{equation} \label{epsilon} 
\epsilon = \epsilon(A/B) =  \sum_{i=1}^{r_1}\dim
  (I_{A/B})_{n_{1,i}} =\sum_{i=1}^{r_1} \bigl[ H_B({n_{1,i}}) -H_A({n_{1,i}})
  \bigr]
\end{equation}
where $H_B$ and $H_A$ are the Hilbert functions of $B$ and $A$. If $B$ is a
complete intersection (CI), allowing $R=B$, then $I_B/I_B^2$ and the normal
module $N_B = (I_B/I_B^2)^*$ are $R$-free of rank $r_1 \geq 0$, and 
\begin{equation*} \label{epsilonCI} 
 \dim_{(B)} \GradAlg^{H_B}(R)=\dim (N_B)_0 = \sum_{i=1}^{r_1} H_B({n_{1,i}})
\end{equation*} 

In Lemma~\ref{psmooth} and Theorem~\ref{mainzero} of the next section we look
upon the special case $B=R$ as a CI with $r_1=0$. Throughout we pass to small
letters to denote the $k$-vector space dimension of the (co)homology groups
involved, e.g. for any $i \geq 0$,
\[
\hho^i(\widetilde{M}) = \dim \HH^i(\widetilde{M}), \ _v\!\hho^i (R,B,M)= \dim
\,_v\!\HH^i (R,B,M), \ _{v}\!\ext_B^i(M,N)= \dim \,_{v}\!\Ext_B^i(M,N).
\]

\begin{lemma} \label{psmooth}
  Let $R \rightarrow B = R/I_B \rightarrow A \simeq B/I_{A/B}$ be
  graded morphisms, let $c=\dim R - \dim B$ and  suppose either
  
  (a) \ \ \ \ $I_B$ is generated by a regular sequence (allowing $R = B$), or
  
  (b) \ \ \ \ $ c \geq 2$ and $B_v \rightarrow A_v$ is an isomorphism for all
  $v \leq
  \max_{i}\{n_{2,i}\}$. \\
  Then $\ _0\!\HH^2 (R,B,I_{A/B}) =0 \ $ and $ \ _0\!\hom_R(I_B,I_{A/B})=
  \epsilon(A/B)$. Moreover $\epsilon(A/B)=0$ if (b) holds.
\end{lemma}
\begin{proof} If $B$ is a CI, then it is well known that $\ _0\!\HH^2
  (R,B,I_{A/B}) =0 $, and moreover that $\ _0\!\HH^1 (R,B,I_{A/B}) \simeq \
  _0\!\Hom_R(I_B/I_B^2,I_{A/B}) \simeq (\oplus_i I_{A/B}(n_{1,i}))_0$ and we
  get the lemma in this case. In (b) it suffices by \eqref{cohom} to show $\
  _0\!\Ext_R^i(I_B,I_{A/B})=0$ for $i \leq 1.$ Applying $\
  _0\!\Hom_R(-,I_{A/B}) $ to the minimal resolution of $I_B$ deduced from
  \eqref{resB}, we conclude by the assumptions of (b).
\end{proof}

The following Proposition 
is a part of Prop.\! 4 of \cite{K04} and is used quite often in this paper.
Below $\GradAlg(H_B,H_A)$ is the representing object of the functor deforming
surjections $B \rightarrow A$ of graded quotients of $R$ of positive depth
(for non-Artinian quotients) and with Hilbert functions $H_B$ and $H_A$ of $B$
and $A$ respectively. Then there exist natural projection morphisms
$p:\GradAlg(H_B,H_A) \rightarrow \GradAlg(H_A)$, induced by $p((B \rightarrow
A))=(A)$, and $q:\GradAlg(H_B,H_A) \rightarrow \GradAlg(H_B)$, induced by
$q((B \rightarrow A))=(B)$, which under the assumptions of Prop.\! 4 of
\cite{K04} have nice properties. Recall that $A$ is called {\it $H_B$-generic}
if there is an {\em open} subset $U \ni (A)$ of $\GradAlg^{H_A}(R)$ such that
every $(A') \in U$ belongs to $\im p$.
Now since the surjectivity of the natural map $\ _0\!\Hom_B(I_{B},B)
\rightarrow \ _0\!\Hom_R(I_{B},A)$ together with the injectivity of $
_0\!\HH^2 (R,B,B) \rightarrow \ _0\!\HH^2 (R,B,A)$ is equivalent to $$\ 
_0\!\HH^2 (R,B,I_{A/B}) =0 \ $$
by the long exact sequence of algebra
cohomology, we may state \cite{K04}, Prop.\! 4 (i), resp. (ii) as (i), resp.
(ii) of the Proposition below.

\begin{proposition}\label{prsmoothness} 
  Let $B = R/I_B \rightarrow A \simeq B/I_{A/B}$ be a graded morphism of
  quotients of $R$. 
  
  (i) If $\ _0\!\HH^2(B,A,A)=0$, (e.g. $ _0\!\Ext_B^1(I_{A/B},A) =0$), then
  the projection $q: \GradAlg(H_B,H_A) \rightarrow \GradAlg(H_B)$ induced by
  $q((B \rightarrow A))=(B)$ is smooth with fiber dimension $\ 
  _0\!\hom_B(I_{A/B},A)$ at $(B \rightarrow A)$.
  
 (ii) If $_0\!\HH^2 (R,B,I_{A/B}) =0$,
 then the projection $p: \GradAlg(H_B,H_A) \rightarrow \GradAlg(H_A)$ induced
 by $p((B \rightarrow A))=(A)$ is smooth with fiber dimension $\ 
 _0\!\hom_R(I_B,I_{A/B})$ at $(B \rightarrow A)$. In particular $A$ is
 $H_B$-generic.
\end{proposition}

\begin{corollary}\label{corprsmooth} 
  Let $B \rightarrow A$ be a graded surjection of quotients of $R$.  If
  $\ _0\!\HH^2(B,A,A)=0$ and $_0\!\HH^2 (R,B,I_{A/B}) =0$, then $$ \ \dim
  (N_A)_0 + \ _0\!\hom_R(I_B,I_{A/B}) = \dim (N_B)_0 + \ _0\!\hom_B(I_{A/B},A)
  \ , \ {\rm and}
  $$
 $$ \
 \dim_{(A)} \GradAlg(H_A) + \ _0\!\hom_R(I_B,I_{A/B}) = \dim_{(B)}
 \GradAlg(H_B)+ \ _0\!\hom_B(I_{A/B},A) \ . \
  $$
  Hence $A$ is unobstructed as a graded $R$-algebra iff $B$ is
  unobstructed as a graded $R$-algebra.
%
\end{corollary}

\begin{proof} Using Proposition~\ref{prsmoothness}(i) we get $\dim_{(B
    \rightarrow A)} \GradAlg(H_B,H_A) = \dim_{(B)} \GradAlg(H_B)+ \ 
  _0\!\hom_B(I_{A/B},A)$ while (ii) implies $\dim_{(B \rightarrow A)}
  \GradAlg(H_B,H_A) = \dim_{(A)} \GradAlg(H_A) + _0\!\hom_R(I_B,I_{A/B})$
  which gives one of the dimension formulas. Since smooth morphisms imply
  surjective tangent maps of their tangent spaces and since the $\Hom$-groups
  of Proposition~\ref{prsmoothness} are the tangent spaces of the fibers, we
  can argue as above to get the other dimension formula.
\end{proof}
  
\section{Families of one dimensional $R$-quotients.} 

In this section we focus on families of zero schemes in $\pp = \pp^{n-1}$ with
fixed Hilbert function $H$, i.e. we study the (possibly non-reduced)
postulation Hilbert scheme $\GradAlg^{H}(R)$ where $H(v)$ is a constant for $v
>> 0$. 
 
If $Y \subset \pp=\Proj R$ is a closed subscheme and $X=\Proj(A)$ is obtained
by choosing s general points on $Y=\Proj(B)$ (in the sense of Geramita et al.
\cite{GMR}), the main theorem of this section implies that $A$ and $B$ are
simultaneously unobstructed as graded algebras and $\dim_{(A)}
\GradAlg^{H_A}(R)$ and $\dim_{(B)}\GradAlg^{H_B}(R)$ are closely related
(Theorem~\ref{mainzero}, Corollary~\ref{sgeneric}). Even though this result
may seem new as stated, it is a straightforward consequence of Theorem 9.16 of
\cite{KMMNP} if Y is a curve. In this section we generalize the result to any
scheme $Y$. In Proposition~\ref{sgen} 
we extend the result to families, and we finish by a theorem on
linkage of families.

A zero-dimensional closed scheme $X \hookrightarrow Y$ is said to be a local
complete intersection (l.c.i) of codimension $(r_1,...,r_t)$ with respect to
$X=X_1 \cup ... \cup X_t$ if $X$ can be written as a disjoint union $X=X_1
\cup ... \cup X_t$ where, for each $i$, the ideal $(\sJ_{X/Y})_{,x}$ is
generated by an $\sO_{Y,x}$-regular sequence of length $r_i$ for every $x \in
X_i$. If $r_i$ are equal for all $i$, say $r_i=r$, we simply say $X
\hookrightarrow Y$ is an l.c.i of codimension $r$. Note that in the case
$r_i=0$, then $X_i$ is mapped isomorphically onto an open subscheme of $Y$.
Below $N_B:=\Hom_B(I_B/I_B^2,B)$ is the normal module of $B$ in $R$, and
$\epsilon(A/B)$ is defined in \eqref{epsilon} and $n_{2,j}$ in \eqref{resB}.

\begin{theorem} \label{mainzero} Let $R$ be a polynomial $k$-algebra and let
  $B = R/I_B \rightarrow A = R/I_A$ be a graded morphism such that $A$ is
  Cohen-Macaulay of dimension one and $\depth_{\mathfrak{m}}B \geq 1$, and
  such that $X:= \Proj(A) \hookrightarrow Y:=\Proj(B)$ is a local complete
  intersection of codimension $(r_1,...,r_t)$ with respect to $X=X_1 \cup ...
  \cup X_t$. Let $H_A(v)= s$ and $H_{X_i}(v)= s_i$ for $v >>0$ (so $s=
  \sum_is_i $) and suppose either
  
  (a) \ \ \ \ $I_B$ is generated by a regular sequence (allowing $R = B$), or
  
  (b) \ \ \ \ $B_v \rightarrow A_v$ is an isomorphism for all $v \leq
  \max_{i}\{n_{2,i}\}$ and $\dim R - \dim B \geq 2$. \\
  Moreover suppose there is an integer $j$ such that $B_v \simeq A_v$ for all
  $v \leq j-1$ and such that $I_A$ is $(j+1)$-regular (i.e. $\reg(A) \leq j$).
  Then $ \ \dim (N_A)_0= \dim (N_B)_0 + \sum_ir_is_i - \epsilon(A/B) \ , $ and
  $$ \ \dim_{(A)} \GradAlg^{H_A}(R)= \dim_{(B)} \GradAlg^{H_B}(R)+
  \sum_i^tr_is_i- \epsilon(A/B) \ .$$ In particular $A$ is unobstructed as a
  graded $R$-algebra if and only if $B$ is unobstructed as a graded
  $R$-algebra.
\end{theorem}

\begin{remark} \label{linres} Theorem~\ref{mainzero} applies to quotients $B
  \rightarrow A \simeq B/I_{A/B}$ where the mapping cone construction produces
  the minimal resolution of A from the free resolution of $B$ and a {\rm
    semi-linear} resolution of $I_{A/B}$ (modulo redundant terms). For
  instance if $M:=I_{A/B}$ and $t=r-1$ in \eqref{linear} (i.e. $ \depth
  I_{A/B} =2$), then $B_v \simeq A_v$ for all $v \leq j-1$ and $I_A$ is
  $(j+1)$-regular, and Theorem~\ref{mainzero} applies. Also in the case $
  \depth I_{A/B} =1$ in which the ``${G}_{r-1} \rightarrow . . .  \rightarrow
  {G}_1$''-part of the minimal resolution $0 \rightarrow {G}_{r} \rightarrow .
  . . \rightarrow {G}_1 \rightarrow I_{A/B} \rightarrow 0$ is semi-linear,
  then Theorem~\ref{mainzero} applies because the contribution from
  $G_{r}$ becomes redundant in the minimal resolution of $A$. Thus the
  condition on $j$ of Theorem~\ref{mainzero} essentially impose on $I_{A/B} $
  to have a semi-linear resolution, which, in the non CI case, must be large
  enough to have (b) fulfilled (e.g. $j \geq \reg(I_B)+2$).
\end{remark}
\begin{proof} It is enough to prove the two dimension formulas. Due to
  Corollary~\ref{corprsmooth} it suffices to show $\ _0\!\HH^2(B,A,A)=0$ and $
  _0\!\hom_B(I_{A/B},A) = \sum_ir_is_i $, as well as $\ _0\!\HH^2
  (R,B,I_{A/B})=0$ and $ \ _0\!\hom_R(I_B, I_{A/B}) = \epsilon(A/B) $. The
  latter follows from Lemma~\ref{psmooth}. Let $\sN_{X/Y}$ the normal sheaf of
  $X \hookrightarrow Y$. Since $\dim X = 0$ and the composition $x
  \hookrightarrow X \hookrightarrow Y$ is a local complete intersection for
  any $x \in X$, then the exact sequence in the proof of Thm.  9.16 of
  \cite{KMMNP} shows $\ _0\!\HH^2(B,A,A)=0$ and $ _0\!\Hom_B(I_{A/B},A) \simeq
  \ _0\!\HH^1(B,A,A) \simeq \HH^0(\sN_{X/Y})$ provided $$
  \ 
  _0\!\Hom_R(I_{A/B},\HH_{\mathfrak m}^1(A) ) = 0\ . $$
  Since $\ \ 
  \HH_{\mathfrak m}^1(A)_v \simeq \ \HH_{\mathfrak m}^2(I_A)_v \simeq \ 
  _{-v}\!\Ext_R^{n-2}(I_A,R(-n))^{\vee}$, we get that $\ \HH_{\mathfrak
    m}^1(A)_v = 0$ for $v \geq j$ by the $(j+1)$-regularity of $I_A$. Using
  $(I_{A/B})_{j-1}=0$ we conclude that $ _0\!\Hom_R(I_{A/B},\HH_{\mathfrak
    m}^1(A) ) = 0. $
  
  Hence it suffices to show $\ \dim \HH^0(\sN_{X/Y})= \sum_ir_is_i$.
  Now since $\Supp(X)$ is finite, we know that $\hho^0(\sO_{X_i})=
  \sum_{x \in \Supp(X_i)} \length(O_{X_i,x}) =s_i$. Using that
  $\sN_{X/Y,x}$ is a free $\sO_{X,x}$-module of rank $r_i$ for any $x
  \in \Supp(X_i)$, we conclude by
  $$\hho^0(\sN_{X/Y})= \sum_i^t  \sum_{x \in \Supp(X_i)}
  \length(\sN_{X/Y,x}) = \sum_i^t \sum_{x
    \in \Supp(X_i)} r_i \cdot \length(O_{X_i,x})= \sum_i^t r_is_i \ .$$
\end{proof}

Moreover Theorem~\ref{mainzero} is precisely what we need to treat the
case where $X$ consists of $s$ (distinct) points in generic position
on $Y$ (i.e. $H_A$ is the truncation of $H_B$ at the level $s$).
Indeed Geramita-Maroscia-Roberts (\cite{GMR}) defines such a
truncated Hilbert function by
$$H_A(i) = \inf \{H_B(i), s\} \ , \ $$
and they show that there exists
a reduced scheme $X$ on $Y$ with truncated Hilbert function $H_A$
provided $Y$ is reduced and consists of more that $s$ points. We get

 \begin{corollary} \label{sgeneric} Let $Y=\Proj(B)$, $B=R/I_B$, be a reduced
   scheme consisting of more than $s$ points, and let $X=\Proj(A)$ be $s$
   points (avoid $SingY$) in generic position on $Y$.
   Let $j$ be the smallest number such that $H_A(j) \neq H_B(j)$. If $Y$ is
   not a CI, suppose $j \geq \max_{i}\{n_{2,i}\} +1$ (e.g. $j \geq
   \reg(I_B)+2$).  Then $X \hookrightarrow Y$ is an l.c.i of codimension
   $(r_1,...,r_t)$ with respect to some decomposition $X=X_1 \cup ... \cup
   X_t$. Moreover $ \ \dim (N_A)_0= \dim (N_B)_0 + \sum_i r_is_i \ -
   \epsilon(A/B) \ , $ and $$
   \ \dim_{(A)} \GradAlg^{H_A}(R)= \dim_{(B)}
   \GradAlg^{H_B}(R)+ \sum_i^t r_is_i \ - \epsilon(A/B) \ .$$
   Hence $A$ is
   unobstructed as a graded $R$-algebra iff $B$ is unobstructed as a graded
   $R$-algebra.
\end{corollary}

\begin{proof}
  Since $\sO_{Y,x}$ and $\sO_{X,x}$ are regular local rings for any $x
  \in X$, it follows that $X \hookrightarrow Y$ is an l.c.i of
  codimension as in Corollary~\ref{sgeneric}.  By the definition of $s$
  generic points, $(I_{A/B})_{j-1}=0$. Since $H_A(v) \neq H_A(j)$ for
  $v\geq j$, it follows that $I_A$ is $(j+1)$-regular (and $j$-regular
  if $s =H_A(j-1)$). Then Theorem~\ref{mainzero} applies supposing $j$
  large enough.
\end{proof}

\begin{remark} \label{nnn} It is well known that the Hilbert polynomial
  $p_B(x)$ equals $H_B(x)$ for all $x \geq \reg(I_B)-1$. Thus the number $j
  \geq \reg(I_B)+2$ of Corollary~\ref{sgeneric} is so large that
  $p_B(x)=H_B(x)$ for $x \geq j-3$. In particular using
  Corollary~\ref{sgeneric} with say $j = \reg(I_B)+2$, we get an algebra $A$
  with Hilbert function $H_A(x)=H_B(x) = p_B(x)$ for $x \in \{j-3,j-2,j-1\}$
  and $H_A(x)=s$ for $x \geq j$.
\end{remark}
\begin{example} \label{exsing} (an obstructed one-dimensional level 
  algebra with $h$-vector $(1,4,4,4,0,0,...)$). 
  
  The part of the Hilbert scheme $H:=\Hilb^{dx+1-g}(\pp^4)$ consisting of
  rational normal curves of degree $d=4$ is thoroughly studied (\cite{MDPi},
  \cite{NS}). It is generically smooth and its closure forms an irreducible
  component $V$ of $H$ of dimension $5d+1 = 21$. The normal sheaf of the
  general curve $Y_g$ satisfy $ \HH^1(\sN_{Y_g})=0$, while for instance
  $Y=\Proj(B)$; the union of four lines meeting at a point, belongs to the
  same component $V$ and satisfies $\dim \HH^1(\sN_{Y})=3$ (cf. \cite{KMMNP},
  Rem.\! 9.9), i.e.  $Y$ is an obstructed reduced arithmetically CM (ACM)
  curve. Both curves have the same graded Betti numbers, e.g.
\begin{equation*} \label{resolB}
   0 \rightarrow R(-4)^3 \rightarrow R(-3)^8 \rightarrow R(-2)^6 \rightarrow
   R \rightarrow B \rightarrow 0 \ .
\end{equation*}
Since the locus of ACM curves in $ \GradAlg(H)$ is open in $H$ by
Proposition~\ref{Grad}, then $V$ corresponds to an irreducible component of $
\GradAlg(H)$ to which $(Y_g)$ and $(Y)$ belong. Let $X= \Proj(A)$ (resp.
$X_g=\Proj(A_g)$) be obtained by choosing $s \geq 13$ generic points on $Y$
(resp. $Y_g$). Since $\dim B_v =4v+1$ for $v \geq 0$, we see that
Corollary~\ref{sgeneric} applies for $j \geq 4$. It follows that $A_g$ is
unobstructed while $A$ is obstructed as graded $R$-algebras and
$$\dim_{(A_g)} \GradAlg(H_{A_g}) = \dim (N_{A_g})_0 = \hho^0(
{\sN}_{Y_g})+ s =21+s$$
(resp. $ \ \dim (N_A)_0= \dim (N_B)_0 + s
=24+s$). If $s=13$ it is straightforward to see that $A_g$ and $A$ are
level algebras with the same graded Betti numbers, e.g. its minimal
resolution is
   \begin{equation*}
   \small
   0 \rightarrow R(-7)^4 \rightarrow R(-6)^{12} \oplus R(-4)^3 \rightarrow
   R(-5)^{12}\oplus R(-3)^8 \rightarrow R(-4)^4 \oplus R(-2)^6 \rightarrow
   I_A \rightarrow 0.
   \end{equation*}
\end{example}

Corollary~\ref{sgeneric} applies also to families of reduced schemes $Y$ which
is not necessarily equidimensional.

 \begin{example} \label{nonequi} Let $H(x) = 3x+1$ for $x \geq 0$, so
   $H=(1,4,7,10,13,...)$. If $Y_1 = \Proj(B_1) \subset \pp^3$ is a twisted
   cubic curve and $Y_2 = \Proj(B_2)$ is the union of a plane space curve $C$
   of degree $3$ and a point $P$ outside the plane containing $C$, then it is
   easy to see that both curves belong to the same stratum $\GradAlg(H)$
   of the Hilbert scheme $\Hilb^{3x+1}(\pp^3)$. We claim they belong to two
   different components of $\GradAlg(H)$. Indeed $(Y_1)$ belongs to a
   12-dimensional irreducible component of $\GradAlg(H)$, and using $
   \sN_{Y_2} \simeq \sN_C \oplus \sN_P$ and that $C \hookrightarrow \pp^3$ and
   $P \hookrightarrow \pp^3$ are CI, we easily get $ \hho^0(\sN_{Y_2}) =15$
   and $ \HH^1(\sN_{Y_2})=0$. Invoking Proposition~\ref{Grad} we see that
   $(Y_2)$ belongs to a 15-dimensional irreducible component of
   $\GradAlg(H)$, cf. \cite{PiS} for a complete description of $
   \Hilb^{3x+1}(\pp^3)$. The minimal resolution of $I_{B_2}$ (resp. $I_{B_1}$)
   is of the form
\begin{equation} \label{resolutB}
   0 \rightarrow R(-4) \rightarrow R(-4) \oplus R(-3)^3 \rightarrow
   R(-3) \oplus R(-2)^3
   \rightarrow I_{B_2} \rightarrow 0 \ 
\end{equation} 
(resp. of the form \eqref{resolutB} where both $R(-4)$ and two of
$R(-3)$ are removed). Since Corollary~\ref{sgeneric} applies for $j
\geq 5$, let $X_1= \Proj(A_1)$ (resp. $X_2=\Proj(A_2)$) be obtained by
choosing $s \geq 13$ generic points on $Y_1$ (resp. $Y_2$); on $Y_2$
we must choose $P$ as one of the $s$ generic points to get the right
Hilbert function. It follows that $A_i$ are unobstructed as graded
$R$-algebras for $i=1$ and $2$ and that $\dim_{(A_1)} \GradAlg(H') =
12+s$ where $H'=(1,4,7,...,3j-2,s,s,...)$, $3j-2 \leq s < 3j+1$. Since
$X_2 \hookrightarrow Y_2$ is an l.c.i of codimension $(1,0)$ with
respect to the decomposition $X_2 = C_2 \cup P$ where $C_2$ consists
of $s-1$ points, we get
$$ \
\dim_{(A_2)} \GradAlg(H')= 15 + s-1 =14+s \ .$$ Hence we get two different
components of $ \GradAlg(H')$. Finally if $s=13$ it is straightforward to see
that $A_2$ have the minimal resolution
   \begin{equation*}
   \small
   0 \rightarrow R(-7)^{3} \oplus R(-4) \rightarrow
   R(-6)^{6} \oplus R(-4)\oplus R(-3)^3 \rightarrow R(-5)^3 \oplus
   R(-3)\oplus R(-2)^3  \rightarrow  I_{A_2} \rightarrow 0.
   \end{equation*}
\end{example}

Once the connection between $\GradAlg(H_A)$ and $
\GradAlg(H_B)$ for $s$ generic points $X$ on $Y$ is that nice as
described in Corollary~\ref{sgeneric}, one may also ask if their
irreducible components correspond exactly and similar questions. E.g.,
may we look upon $A_g$ of Example~\ref{exsing} as the general element
of an irreducible component of $\GradAlg(H_{A_g})$? To see the
answer is yes we use some ideas of \cite{K04}.

\begin{definition}
  Inside $\GradAlg(H)$ we look to the following open subsets, $\Smc(H)$ (resp.
  $\SmCM(H)$), consisting of points $(R/I)$ such that $\Proj(R/I)$ is a smooth
  geometrically connected scheme (resp. smooth and arithmetically CM). Here
  ``points'' should be considered as ``$\Omega$-points'' where $\Omega$ is an
  overfield of $k$.  Moreover let $\Smc(H)_n$ be the open subset of $\Smc(H)$
  consisting of points $(R/I)$ where the Castelnuovo-Mumford regularity
  satisfies $\reg(I) \leq n$.  Similarly we let $\CI(H)$ (resp. $\CM(H)$)
  consist of points $(R/I)$ where $I$ is generated by a regular sequence
  (resp. $R/I$ is CM).
\end{definition}
Now let $$ \SmCM(H_B,H_A)_n := p^{-1}(\SmCM(H_A)) \cap q^{-1}(\Smc(H_B)_n)$$
where $ q: \GradAlg(H_B,H_A) \rightarrow \GradAlg(H_B) $ and $p:
\GradAlg(H_B,H_A) \rightarrow \GradAlg(H_A)$ are the two natural projection
morphism. (e.g. $q((B \rightarrow A))= (B)$). Denoting the following
restrictions of $p$ and $q$ by the same letters, we get a diagram (incidence
correspondence)
\begin{equation} \label{incid}
  \begin{array}[h]{ccc}
    \SmCM(H_B,H_A)_n & \stackrel{q}{\longrightarrow} & \Smc(H_B)_n \subset
    \GradAlg(H_B) \\
    \rule[-2mm]{0pt}{6mm} \downarrow ^{p} \\
    \GradAlg(H_A)
  \end{array}
\end{equation}

\begin{proposition} \label{sgen} Let $H_B$ be the Hilbert function of some
  smooth connected curve and let $H_A$ be its truncated Hilbert
  function at the level $s$. Let $ j = \min\{ i \arrowvert H_A(i) \neq
  H_B(i) \} , $ let $n \leq j-2$ and look to the maps $p$ and $q$ in
  \eqref{incid}. Then

  (i) \ \ \  $q$ is smooth and surjective and its fibers are geometrically
  connected, of fiber dimension $s$, and

  (ii) \ \ \ $p$ is an isomorphism onto an open subscheme of $\
  \GradAlg(H_A)$. \\ 
  In particular the correspondence \eqref{incid} determines a well-defined
  injective application $\pi$ from the set of irreducible components $W$ of
  $\ \Smc(H_B)_n$, to the set of irreducible components $V$ of $\
  \GradAlg(H_A)$, in which generically smooth components correspond.
  Indeed $V=\pi(W)$ is the closure of $p(q^{-1}(W))$, and we have
  $$
  \ \dim V = \ \dim W + s \ .$$
\end{proposition}

\begin{proof}
  (i) By Geramita et al. \cite{GMR} we get the surjectivity of $q$.
  Since we showed $\ _0\!\HH^2(B,A,A)=0$ in Theorem~\ref{mainzero},
  the smoothness of $q$ follows immediately from
  Proposition~\ref{prsmoothness}(i). To show that the fibers of $q$
  are (geometrically) connected, one may simply look to the fiber as
  the variation of $s$ generic points on a fixed $Y$, i.e. as a
  non-empty dense set of $Y^s$. This set is irreducible since Y is
  irreducible, and we conclude as claimed. 
   
  (ii) In Lemma~\ref{psmooth} we showed $\ 
  _0\!\Ext_R^i(I_B,I_{A/B})=0$ for $i \leq 1$ assuming $j \geq
  \reg(I_B)+2$. By Proposition~\ref{prsmoothness}(ii) this implies
  that $p$ is smooth and unramified. It is easy to see that $j \geq
  \reg(I_B)+1$ implies that $p$ is injective (in fact, universally
  injective or ``radiciel''), cf. Lemma 7(a) of \cite{K04}. Hence we
  get (ii) by \cite{EGA}, Thm.\!\ 17.9.1.  Now combining (i) and
  \cite{HAR}, Prop.\! 1.8, we get that $q^{-1}(W)$ is an irreducible
  component of $ \SmCM(H_B,H_A)_n$. The application $\pi$ is therefore
  well defined, and it is injective by (ii). Finally since $q$ is
  smooth and $p$ is an open immersion, we easily get the dimension
  formulas.
\end{proof}

\begin{remark}
  If we in Proposition~\ref{sgen} drop the assumption $\dim
  \Proj(B)=1$ and maintain the other assumptions, we still get that
  $q$ is smooth and that $p$ is an isomorphism onto an open subscheme
  (but the irreducibility of $q^{-1}(W)$ may fail).
\end{remark}

 
Now we consider an example of several components of $\ \GradAlg(H_A)$, which
one may, as in \cite{Mi}, distinguish by the incomparability of the set of
graded Betti numbers (Remark~\ref{conn}). Applying, however,
Proposition~\ref{sgen} to our example we can also describe well the graded
Betti numbers of some algebras in the {\it intersection of the two components}.

 \begin{example} \label{twocomp} In \cite{W1} C. Walter gives examples
   of infinitely many Hilbert schemes of space curves containing {\rm
     obstructed} smooth curves of maximal rank. Indeed his example of
   a smooth space curve $Y$ of the lowest degree (i.e. the curve with
   Hilbert polynomial $p(x)=33x-116$ which we consider below) was
   independently discovered by Bolondi et al \cite{BKM} and it was the
   first example of an obstructed curve of maximal rank which was
   detected. In \cite{BKM} we showed that $\Hilb^{33x-116}(\pp^3)$
   contains at least two irreducible components whose intersection
   contains $(Y)$. Since the curve $Y=\Proj(B)$ is of maximal rank, we
   have $\ _0\!\Hom_R (I_B,\HH_{\mathfrak m}^1(B)) = 0$, and
   Proposition~\ref{Grad} applies. It follows that the corresponding
   algebra $B$ is {\rm obstructed} as a graded algebra and that $(B)$
   sits in the intersection of two irreducible components $W_1$ and
   $W_2$, both of dimension $4d=132$, of the postulation Hilbert
   scheme of space curves $\ \GradAlg(H_B)$, cf.  \cite{K05}, ex. 35.
   
   In \cite{K05} we also considered the  minimal resolution of $B$ as
   well as the minimal resolution of the general elements $B_1$ and
   $B_2$ of $W_1$ and $W_2$ respectively. Indeed
   $$0 \rightarrow G_3= R(-9) \rightarrow R(-10)^2 \oplus R(-9) \oplus R(-8)^4
   \rightarrow R(-9) \oplus R(-8) \oplus R(-7)^5 \rightarrow I_B \rightarrow
   0$$
   is exact and we get the minimal resolution of $B_1$ (resp.  $B_2$) by
   making the factor $R(-9)$ redundant in two different ways, i.e.  by
   removing this factor from the leftmost ($G_3$) and the middle term ($G_2$),
   making $B_1$ CM, (resp. from $G_2$ and rightmost term $G_1$). The
   Castelnuovo-Mumford regularity for all three curves satisfy $\reg(I)=9$,
   and the Hilbert function of all algebras are
   \[ (1,4,10,20,35,56,84,115,148,181,214,...)\ . \]
   Thus taking $s \geq 214$ points $X=\Proj(A)$ on $Y$ in general position and
   correspondingly for the others, then Proposition~\ref{sgen} applies with $j
   \geq 11$. Or more precisely, both $W_1$ and $W_2$ and its intersection
   essentially belong to $\Smc(H_B)_9 \subset \GradAlg(H_B)$, and
   Proposition~\ref{sgen} applies to (every element $\Proj(B')$ of)
   $\Smc(H_B)_9$ and an $s$-dimensional linear space of choices of $s$ generic
   points on $\Proj(B')$. Hence for each $s \geq 214$ it follows that $X$ is
   in the intersection of two irreducible component $V_1$ and $V_2$ of the
   postulation Hilbert scheme $\ \GradAlg(H_A)$ of dimension $\dim V_i =
   132+s$ for $i=1,2$. In the special case $s=214$ we have $\Delta H_A =
   (1,3,6,10,15,21,28,31,33,33,33,0,...)$, and it is not difficult to see that
   the minimal resolution of $I_A$ is
   $$0 \rightarrow R(-13)^{33} \oplus G_3 \rightarrow R(-12)^{66}
   \oplus G_2 \rightarrow R(-11)^{33} \oplus G_1 \rightarrow I_A
   \rightarrow 0.$$
   and that those of the corresponding general element $\Proj(A_i)$ of
   $V_i$ are obtained by removing the free factor $R(-9)$ from $G_3$
   and $G_2$ (resp. from $G_2$ and $G_1$). Looking to the
   corresponding sets of graded Betti numbers of $A_1$ and $A_2$ we
   see they are incomparable. 
 \end{example}
 
 We finish this section by recalling some known results about the postulation
 Hilbert scheme $ \GradAlg^{H}(R)$, consisting of zero-dimensional schemes
 $\Proj(A)$ of degree $s$. Since we have observed that $_0\! \HH^2 (R,A,A)=0$
 implies the smoothness of $ \GradAlg^{H}(R)$ at $(A)$, we remark that the
 smoothness results of Remark~\ref{summarize}(i) (when $A$ is generically a
 CI) and of Remark~\ref{summarize}(ii) also follow from works of Herzog,
 Buchweitz-Ulrich and Huneke (\cite{HER}, \cite{BU83} and \cite{H84}). Now in
 addition to Theorem~\ref{mainzero} and Proposition~\ref{sgen}, we have

\begin{remark} \label{summarize} 
  
  (i) If $\Proj(R)= \pp^2$, then Gotzmann (\cite{Go}) shows that $
  \GradAlg^{H}(R)$ is irreducible and she finds its dimension (\cite{IK},
  Thm.\!  5.21 and Thm.\! 5.51). It is smooth by licciness and say (iii) below
  (or by \cite{Go} provided $ \GradAlg^{H}(R)$ is reduced). As indicated by
  Iarrobino-Kanev (\cite{IK}, Remark to Thm.\! 5.51), the dimension formula
  given in \cite{KM04}, Rem.\! 4.4, holds in this case (\cite{KM04}, Rem.\!
  4.6).
 
  (ii) If $\Proj(R)= \pp^3$, then the the open part of $ \GradAlg^{H}(R)$
  consisting of Gorenstein quotients is irreducible (cf.\! \cite{Di}), of
  known dimension by (\cite{KM98}, Remark to Thm.\! 2.6) and smooth by say
  (iii) below. This dimension formula is included in \cite{K98}, Thm.\! 2.3
  with a proof (which also takes care of the Artinian case). \cite{K98},
  Prop.\!  3.1 contains a second ``dual'' dimension formula for the same
  parameter space.

 (iii) Let $\Proj(R)= \pp^n$ and let $A$ and $A'$ be two graded CM quotients
 algebraically linked by a CI $B$ of type $(a_1,...,a_m)$ with
 resolution \eqref{resB}. By \cite{K98}, Prop. 1.7, then $A$ and $A'$ are
 simultaneously unobstructed as graded algebras, and we have $$ \dim_{(A)}
 \GradAlg^{H_A}(R) - \sum_{i=1}^m H_A(a_i) = \dim_{(A')} \GradAlg^{H_{A'}}(R)
  - \sum_{i=1}^m H_{A'}(a_i) \ .$$
  
  (iv) Let $B = R/I_B$ be a graded, generically Gorenstein CM quotient with
  canonical module $K_B$ and let $A$ be the Gorenstein algebra given by a
  regular section of $\sigma \in (K_B^*)_t$ for some integer $t$, i.e. given
  by a graded exact sequence $ 0 \rightarrow K_B(-t)
  \stackrel{\sigma}{\longrightarrow} B \rightarrow A \rightarrow 0$.
 
 a) If $B$ is licci, then $A$ is unobstructed as a graded $R$-algebra (indeed
 $ \HH^2 (R,A,A)=0$), and,
  $$
  \ \dim_{(A)} \GradAlg^{H_A}(R)= \ \dim_{(B)} \GradAlg^{H_B}(R) + \ \dim
  (K_B^*)_t -1 -\delta(B)_{-t}$$ where $ \delta(B)_{v} =\
  _{v}\hom_B(I_B/I_B^2,K_{B})-\ _{v}\ext_B^1(I_B/I_B^2,K_{B})$. 

  b) If $\ \Proj(B)$ is locally Gorenstein and $t >> 0$, then $A$ and $B$ are
  simultaneously unobstructed as graded algebras, and the dimension formula of
  a) holds (with $ \delta(B)_{-t} =0$).
  
  This theorem is true in arbitrary dimension of $B$. It is proved in
  \cite{K04}, Thm.\! 16 and is a substantial generalization of some of the
  statements of (ii) above because, when we apply it to a CM $B$ of
  codimension two (necessarily licci), we get the dimension formula of (ii) by
  \cite{K04}, Ex.\! 26. The preprint \cite{K04i} contains further
  generalizations of this theorem.
  \end{remark}
  
  By (iii) we see that CI-linkage preserves the smoothness of the parameter
  spaces. Due to \cite{K88}, Prop.\! 3.4 it also preserves the irreducibility
  of the linked family. To define the linked family, let $ \CICM(H_B,H_A)$
  consist of points $(B \rightarrow A)$ such that $B$ is CI and $A$ is CM,
  i.e. $$
  \CICM(H_B,H_A) := p^{-1}(\CM(H_A)) \cap q^{-1}(\CI(H_B))$$
  where $p:
  \GradAlg(H_B,H_A) \rightarrow \GradAlg(H_A)$ (resp. $q$) is the second
  (resp. first) projection morphism (e.g. $q((B \rightarrow A))= (B)$). In the
  case $\dim A = \dim B$ (not necessarily equal to one) and $(B \rightarrow A)
  \in \CICM(H_B,H_A)$, the linked algebra is defined by $A':=B/ \Hom_B(A,B)$.
  This definition extends to families and preserves flatness \cite{K88}.
  Indeed by \cite{K88}, Thm.\!  2.6 there is an isomorphism $\tau$ of schemes
  and obvious second projection morphisms $p$ and $p'$ fitting into
  \begin{equation} \label{incidlink}
    \begin{array}[h]{ccc}
      \tau : \CICM(H_B,H_A) & \stackrel{\simeq}{\longrightarrow} &
      \CICM(H_B,H_{A'}) \\
      \rule[-2mm]{0pt}{6mm} \downarrow ^{p} && \downarrow ^{p'} \\
      \GradAlg(H_A) && \GradAlg(H_{A'})
    \end{array}
  \end{equation}
  where $\tau$ is given by sending $(B_1 \rightarrow A_1)$ to $(B_1
  \rightarrow A_1':=B_1/ \Hom_{B_1}(A_1,B_1))$.

\begin{definition} \label{linkfam} Let the Hilbert polynomials $p_B$ and $p_A$
  (corresponding to $H_B$ and $H_A$ respectively) have the same degree $( \geq
  -1)$ and let $U$ be a locally closed subset of $\ \im p$ in
  \eqref{incidlink}. Then the $H_B$-linked family of $U$ is $$U':=p'(\tau
  (p^{-1}(U))$$
  \end{definition} 
%

  \begin{theorem} \label{mainlink} In \eqref{incidlink} the morphisms
    $p$ and $p'$ are smooth and its fibers are geometrically
    connected, of fiber dimension $\epsilon(A/B)$ at $(B \rightarrow
    A)$ and $\epsilon(A'/B)$ at $(B \rightarrow A')$ respectively. In
    particular the $H_B$-linked family $ U'$ is irreducible (resp.
    open in $\GradAlg(H_{A'})$) if and only if $ U $ is irreducible
    (resp. open in $\GradAlg(H_{A})$).
\end{theorem}

\begin{proof}
  The proof of \cite{K98}, Prop. 1.7 takes care of the smoothness of $p$ and
  $p'$ and their fiber dimension. It remains to prove the connectedness of the
  fibers since the other conclusions then follow easily. The connectedness is,
  however, a straightforward consequence of the proof of Theorem 1.16 of
  \cite{K88} (that part of the proof doesn't require $\deg p > 0$ and it is
  easy to reformulate it for the Artinian case as well), cf. \cite{MP}, Ch.\!
  VII for similar results.
\end{proof}

Of course Theorem~\ref{mainlink} implies Remark~\ref{summarize}(iii)
above. We may use Theorem~\ref{mainlink} to see that many other
properties are by preserved by linkage. Indeed subsets of $
\GradAlg(H)$ for which the members allow the same sequence of
CI-linkages which ends in a CI, is irreducible. 
It does not mean that the subset of $\GradAlg(H)$ of licci quotients is
irreducible, as the following example shows.

\begin{example} \label{ekslicci} We claim  $\GradAlg^H(k[x,y,z,w])$
  with $ \Delta H=(1,3,6,6,3,1)$ contains (at least) two irreducible
  components whose general elements are licci. Of course, the general
  element of one of the components is an arithmetically Gorenstein
  scheme consisting of $20$ points, with minimal resolution
  \begin{equation*}
   \small
   0 \rightarrow R(-8) \rightarrow
   R(-5)^{4} \oplus R(-4) \rightarrow R(-4) \oplus
   R(-3)^4  \rightarrow  R \rightarrow A_1 \rightarrow 0.
   \end{equation*}
   We get $A_1$ by starting with a CI of type $(1,1,2)$ and then
   perform general CI-linkages of type $(2,3,3)$ and $(4,3,3)$. It
   follows that the component is generically smooth of dimension $44$
   by using Remark~\ref{summarize}(iii), or Theorem~\ref{mainlink}, twice.
   
   To get the other component, we start with a point $\Proj(A)$, i.e.
   a CI of type $(1,1,1)$, and we proceed by performing six general CI
   linkages of type $(1,2,3)$, $(2,2,4)$, $(2,3,4)$, $(3,4,4)$,
   $(3,4,5)$, $(3,3,5)$, in this order. (The first five linkages are
   the same as for the level algebra 64] in the appendix C of
   \cite{GHMS}; hence  $\GradAlg(H)$
  with $ \Delta H=(1,3,6,7,6,2)$ also contains two ``licci''
  components.) We get in this way an open subset $U$ of $\GradAlg(H)$
   of algebras $A_2$ with minimal resolution
\begin{equation*}
   \small
   0 \rightarrow R(-8) \oplus R(-7) \oplus R(-6) \rightarrow 
   R(-7) \oplus R(-6)\oplus  R(-5)^{5} \rightarrow R(-5) \oplus
   R(-3)^4  \rightarrow  R \rightarrow A_2 \rightarrow 0.
   \end{equation*} 
   Since the Betti numbers do not coincide with the general element
   $A_1$ of the other component, the closure of $U$ must be a
   generically smooth component of dimension $44$ by 
   Theorem~\ref{mainlink}.  
  
   (This example holds correspondingly for codimension 3 quotients
   with $h$-vector $(1,3,6,6,3,1)$ in a polynomial ring of any
   dimension).
\end{example}

\section{Families of Artinian $R$-quotients (possibly Gorenstein).} 

In this section we look to families of {\rm Artinian} algebras $A$ of Hilbert
function $H=H_A$, i.e. we study the scheme $ \GradAlg(H)$ in the Artinian case
with a special look to level and Gorenstein Artinian quotients. In particular
we give examples of codimension 4 (resp. 3) quotients where $ \GradAlg(H)$ has
at least two components with a Gorenstein (resp. level) algebra belonging to
the intersection of the two components. Moreover we notice that almost all the
results of the preceding section (cf. \! Remark~\ref{summarize}) are known in
the Artinian case, except Theorem~\ref{mainzero}, whose corresponding Artinian
result is the main new Theorem of this section. Of course there is a few
changes to Remark~\ref{summarize}, mostly concerned with references, and we
include some further results. To summarize,

\begin{remark} \label{summarize2} (i) Iarrobino shows that $
  \GradAlg^{H}(k[x,y])$ is irreducible (\cite{I}, Thm.\! 2.9) and he
  finds the dimension (\cite{I}, Thm.\! 2.12 and Thm.\! 3.13). It is
  smooth by licciness (or by \cite{I}, Thm.\! 2.9 provided $
  \GradAlg^{H}(k[x,y])$ is reduced). Also in this case, the dimension
  formula given in \cite{KM04}, Rem.\! 4.4, holds (by the indicated
  argument of \cite{KM04}, Rem.\! 4.6).
 
  (ii) If $R=k[x,y,z]$, then the open part of $ \GradAlg^{H}(R)$ consisting of
  Gorenstein quotients is irreducible (\cite{Di}) and smooth of known
  dimension (\cite{K98}, Thm.\! 2.3). See also \cite{H84}, Cor\! 4.9 for the
  smoothness.

  (iii)  of Remark~\ref{summarize} holds as stated in
  Remark~\ref{summarize}.

  (iv) of Remark~\ref{summarize} holds as well. One may make a little
  progress to (iv,b) by stating it as:
  
  b) If $\ \Proj(B)$ is a locally Gorenstein zero-scheme of degree $s$ and if
  $t \geq 2 \reg(I_B)$, then $A$ and $B$ are simultaneously unobstructed as
  graded algebras, and the dimension formula of (iv, a) holds (with $
  \delta(B)_{-t} =0$ and $ \dim (K_B^*)_t =s$, cf.\! \cite{K04}, Rem.\! 22).
  We may even state it with the ideas of Proposition~\ref{sgen} as done in
  Theorem~\ref{maintrans} below (\cite{K04}, Prop.\! 23, cf.\! \cite{K04i}
  for further generalizations).
   
  (v) One may, via the Macaulay correspondence, consider the set $PS(s,j,n)$
  of Gorenstein quotients obtained from the set of homogeneous polynomials of
  degree $j$ in the ``dual'' polynomial ring, of the form $$f
  =L_1^j+...+L_s^j \ ,$$
  where $L_i$ are linear forms and $s$ is fixed. If
  $H_A$ (which we denote by $H(s,j,n)$) contains a subsequence of the form
  $(s,s,s)$, then the closure of $PS(s,j,n)$ is a generically smooth
  irreducible component of $ \GradAlg^{H_A}(R)$ of known dimension (\cite{IK},
  Thm.\! 4.10A and Thm.\!  1.61, see Thm.\! 4.13 for similar results when
  $H_A$ does not contain such a subsequence).

  (vi) In the interesting Gorenstein Artinian codimension 4 case, there is a
  structure theorem when $H_A=(1,4,7,h,i,...)$ with $3h-i-17 \geq 0$, allowing
  us to describe well the corresponding (generically smooth) irreducible
  component of $ \GradAlg^{H_A}(R)$ (\cite{IK}). In \cite{JK} Johannes Kleppe
  comes up with classes of generically smooth components of known dimension of
  a similar nature.

  (v) Compressed Artinian algebras of fixed socle degrees belong to an
  irreducible generically smooth component of known dimension by \cite{I84},
  Thm.\! IIB.
  \end{remark}
  
  To accomplish Remark~\ref{summarize2}(iv,b), let $U_t \subset \GradAlg(H')$
  be an open subscheme consisting of points $(B)$ such that $B$ is CM and such
  that $ \Proj(B)$ is a locally Gorenstein zero-scheme of degree $s$
  satisfying $H_B=H'$ and $ \reg(I_B) \leq t/2$. Recall that a regular section
  of $\sigma \in (K_B^*)_t$ defines a graded Gorenstein algebra $A$ given by
  the exact sequence $ 0 \rightarrow K_B(-t)
  \stackrel{\sigma}{\longrightarrow} B \rightarrow A \rightarrow 0$. Let $q:
  \GradAlg(H_B,H_A) \rightarrow \GradAlg(H_B)$ be the first projection and let
  $q^{-1}(U_t)_{reg}$ be the intersection of $q^{-1}(U_t)$ by the space of
  those quotients $(B \rightarrow A)$ which correspond to regular sections of
  $(K_B^*)_t$. Then we have a diagram where we have restricted the two natural
  projection morphisms $q$ and $p:\GradAlg(H_B,H_A) \rightarrow \GradAlg(H_A)$
  to $\ q^{-1}(U_t)_{reg}$:
  \begin{equation} \label{incidgor}
    \begin{array}[h]{ccc}
      q^{-1}(U_t)_{reg} & \xrightarrow{q_{res}} & \ U_t \subset
      \GradAlg(H_B) \\
      \rule[-2mm]{0pt}{6mm} \downarrow ^{p_{res}} \\
      \GradAlg(H_A)
    \end{array}
  \end{equation}

\begin{theorem} \label{maintrans}
  With notations as above, then
  
  (i) \ \ \  $q_{res} : q^{-1}(U_t)_{reg} \rightarrow \ U_t$ is smooth and
  surjective, and its fibers are geometrically connected of fiber
  dimension $s -1$, and
  
  (ii) \ \ \ $p_{res}$ is an
  isomorphism onto an open subscheme of $\ \GradAlg(H_A)$. \\
  In particular the correspondence \eqref{incidgor} determines a
  well-defined injective application $\pi$ from the set of irreducible
  components \ $W$ of $\ U_t$, to the set of irreducible components \ $V$
  of $\ \GradAlg(H_A)$, in which generically smooth components correspond.
  Indeed \ $V=\pi(W)$ is the closure of $p_{res}(q_{res}^{-1}(W))$, and we have
  $$
  \ \dim V = \ \dim W + s-1 \ .$$
\end{theorem} 

\begin{proof}[Proof \normalfont (also of Remark~\ref{summarize2}(iv,b))] These
  results are almost exactly Thm.\! 16 (cf.\! Rem.\! 22) and Thm.\! 24 (cf.\!
  Prop.\!  23) of \cite{K04} with a slight improvement. Indeed in replacing
  ``$t >> 0$'' by ``$t \geq 2 \reg(I_B)$'' we assumed $B$ to be generically
  syzygetic (e.g.  $\Proj(B)$ locally licci) in \cite{K04}, Rem.\! 22 and
  Prop.\!  23(iii) (to see that $p_{res}$ is smooth). Since, however, $K_B(-t)
  \simeq I_{A/B}$ and the $R$-dual of \eqref{resB} is a free resolution of
  $K_B(1)$ ($B$ is CM and one-dimensional), we have $(I_{A/B})_v=0$ for $v
  \geq t+1- \reg(I_B)$. Hence we may use Lemma~\ref{psmooth} and
  Proposition~\ref{prsmoothness}(ii) of this paper to see that $p_{res}$ is
  smooth under the assumption $t \geq 2 \reg(I_B)$ (quite similar to what we
  showed in \cite{K04}, Rem.\! 14). In Theorem~\ref{maintrans} and
  Remark~\ref{summarize2}(iv,b) it suffices therefore to suppose $t \geq 2
  \reg(I_B)$ without requiring $B$ to be generically syzygetic. Note also that
  $q_{res}$ is surjective by \cite{B992}, Thm.\! 3.2, (cf. \cite{K04}, Rem.\!
  22), and it follows that the correspondence \eqref{incidgor} has the stated
  properties.
 \end{proof}  
 
 Now we illustrate Theorem~\ref{maintrans}. The benefit of using
 Theorem~\ref{maintrans} instead of Remark~\ref{summarize2}(iv,b) is clear
 because it is a statement about the whole subscheme $U_t \subset
 \GradAlg(H')$ and not only about a point in $U_t$. E.g. note that if we apply
 (iv,b) to the two components of $\GradAlg(H')$ of Example~\ref{nonequi}, say
 with $s=13$ and $t \geq 10$ to simplify, we get two components of $\PGor(H)$,
 or of $ \GradAlg(H)$, with
  $$H=(1,4,7,10,13,13,...,13,10,7,4,1)$$
  of dimension $37$ and $39$ where the number $13$ occurs $t-7$ times. Such
  components are now well known (\cite{B99}, see also \cite{IS}). Since
  $\GradAlg(H)$ is connected there are graded Artinian quotients belonging to
  the intersection of the components of $\GradAlg(H)$. But are there
  Gorenstein quotients in this intersection? The answer would have been yes if
  the intersection of the two components of $\GradAlg(H')$ of
  Example~\ref{nonequi} contains points $(B)$ such that $ \Proj(B)$ is a
  locally Gorenstein zero-scheme because we then could apply
  Theorem~\ref{maintrans}! We doubt that there exists such a quotient $B$,
  i.e. we expect that the intersection of the two mentioned components of
  $\PGor(H)$ is empty (cf. Piene-Schlessinger's characterization of the
  intersection of the two components described in Example~\ref{nonequid}).
  Here is an example where we somehow control the intersection.
 
\begin{example} \label{twocompgor} (Two components of  $\PGor(H)$ with
  non-empty intersection)
  
  In Example~\ref{twocomp} we showed the existence of an algebra,
  which we now call $B$ whose corresponding point $(B)$ of the
  postulation Hilbert scheme, $\ \GradAlg(H_B)$, sat in the
  intersection of two irreducible components $V_1$ and $V_2$ of $\ 
  \GradAlg(H_B)$ of dimension $\dim V_i = 132+s=346$ for $i=1,2$. The
  element $(B)$ as well as the two general elements $(B_i)$ of $V_i$
  were obtained by taking $s=214$ generic points on certain curves of
  $\Hilb^{33x-116}(\pp^3)$, one of which with minimal resolution
  $$0 \rightarrow G_3= R(-9) \rightarrow G_2=R(-10)^2 \oplus R(-9)
  \oplus R(-8)^4 \rightarrow G_1=R(-9) \oplus R(-8) \oplus R(-7)^5
  \rightarrow I$$
  Moreover $H_B =
  H_{B_i}=(1,4,10,20,35,56,84,115,148,181,214,214,214,...)$ and the
  minimal resolution of $I_B$ (resp. of $I_{B_1}$, or $I_{B_2}$) was
\begin{equation} \label{G33}  
0 \rightarrow R(-13)^{33} \oplus G_3 \rightarrow R(-12)^{66}
  \oplus G_2 \rightarrow R(-11)^{33} \oplus G_1 \rightarrow I_B
  \rightarrow 0
\end{equation}
  (resp. \eqref{G33} in which the factor $R(-9)$
  from $G_3$ and $G_2$, or from $G_2$ and $G_1$, were removed).
   
  Since we have $\reg(I_B) = \reg(I_{B_i}) = 11$, we may use
  Theorem~\ref{maintrans} to get, for every $t \geq 22$, two generically
  smooth irreducible components of $\PGor(H_A)$ of dimension $132+s+s-1=559$
  whose intersection is non-empty, i.e. the intersection contains an {\rm
    obstructed} Gorenstein Artinian algebra whose $h$-vector is the
  $(t+1)$-tuple $$H_A
  =(1,4,10,20,35,56,84,115,148,181,214,214,...,214,181,148,...,4,1)$$
  where
  the number $214$ occurs $t-19$ times. The corresponding sets of graded Betti
  numbers of the general elements, $A_1$ and $A_2$, of the two components turn
  out to be incomparable because the factors $R(-9)$ (and $R(-t+5)$) appearing
  in the resolution of $I_A$ becomes redundant in different ways in the
  resolutions $I_{A_1}$ and $I_{A_2}$. Of course, for every $s \geq 214$ we
  can construct similar examples.
 \end{example}

 Now we prove the analogue of Theorem~\ref{mainzero} which is the main
 result of this section.

\begin{theorem} \label{mainartin} Let $R$ be a polynomial $k$-algebra and let
  $B = R/I_B \rightarrow A = R/I_A$ be a graded morphism
  such that $A$ is Artinian and $\ \depth_{\mathfrak{m}}B \geq \min(1,\dim
  B)$, and suppose either
  
  (a) \ \ \ \ $I_B$ is generated by a regular sequence (allowing $R = B$), or
  
  (b) \ \ \ \ $B_v \rightarrow A_v$ is an isomorphism for all $v \leq
  \max_{i}\{n_{2,i}\}$  and $\dim R - \dim B \geq 2$. \\ 
  Let $F$ be a free $B$-module such that $F \rightarrow I_{A/B}$ is surjective
  and minimal, and suppose there is an integer $j$ such that the degrees of
  minimal generators of the $B$-module $\ker(F \rightarrow I_{A/B})$ $> j$
  (e.g. $B_v \simeq A_v$ for all $v \leq j-1$) and such that $I_A$ is
  $(j+1)$-regular (i.e.  $A_{j+1} = 0$). Then $ \ \dim (N_A)_0= \dim (N_B)_0 +
  \ _0\hom_B(F,A) - \epsilon(A/B) \ , $ and
 $$
  \ \dim_{(A)} \GradAlg^{H_A}(R)= \dim_{(B)} \GradAlg^{H_B}(R)+ \
  _0\! \hom_B(F,A) - \epsilon(A/B) \ .$$
  In particular $A$ is
  unobstructed as a graded $R$-algebra if and only if $B$ is
  unobstructed as a graded $R$-algebra.
\end{theorem}

At least in the case $\ \depth_{\mathfrak{m}}B \geq 1$, the natural
application of Theorem~\ref{mainartin} is the same as for
Theorem~\ref{mainzero}; the minimal resolution of $A$ should be the one of $B$
in addition to a semi-linear contribution coming from $I_{A/B}$ via the
mapping cone construction, cf. Remark~\ref{lires}.
\begin{proof}
  All we need to prove is the two dimension formulas. Due to
  Corollary~\ref{corprsmooth} it suffices to show $\ 
  _0\!\HH^2 (R,B,I_{A/B})=0$ and $ _0\!\hom_R(I_B,I_{A/B})=
  \epsilon(A/B)$ together with 
\begin{equation} \label{HBA}
  \dim\ _0\! \Hom_B(I_{A/B},A)=\ _0\!
  \hom_B(F,A) \ \ {\rm and} \ \ \ _{0}\!\Ext_B^1(I_{A/B},A)=0 \ , 
\end{equation}
because the latter of \eqref{HBA} implies $\ _0\!\HH^2(B,A,A)=0$. By
Lemma~\ref{psmooth} it suffices to prove \eqref{HBA}. Let
\begin{equation} \label{Bres}
  F' \rightarrow F \rightarrow I_{A/B} \rightarrow 0
\end{equation} 
be the first terms of a $B$-free minimal resolution of $I_{A/B}$.
Applying $ \ _{0}\!\Hom_B(-,A)=0$ onto \eqref{Bres} and using
$A_{j+1}=0$, we get $ \ 
_{0}\!\Hom_B(I_{A/B},A) \simeq \ _{0}\!\Hom_B(F,A)$ and $ \ 
_{0}\!\Ext_B^1(I_{A/B},A)=0$ by the assumption, and we are done.
\end{proof}

\begin{remark} \label{H2=0} By the long exact sequence of algebra cohomology,
  we have the exact sequence
 $$\rightarrow \ _0\!\HH^2(B,A,A)  \rightarrow \ _0\!\HH^2(R,A,A)
 \rightarrow \ _0\!\HH^2(R,B,A) \rightarrow \ . $$ Since it is well known that
 $ \HH^2(R,B,-)=0$ if Theorem~\ref{mainartin}(a) holds and since we have $\
 _0\!\HH^2(B,A,A)=0$ by the proof above, it follows that we in
 Theorem~\ref{mainartin}(a) have $$\ _0\!\HH^2(R,A,A)=0 \ . $$
 \end{remark}

 \begin{remark} \label{reg+2} A natural choice of $j$ in
   Theorem~\ref{mainartin} such that $(I_{A/B})_{j-1}=0$ and such that (b)
   holds, is $j \geq \reg(I_B)+2$, in which case we have that $I_A$ is
   $(j+1)$-regular iff $I_{A/B}$ is $(j+1)$-regular, and that $H_A(x)=H_B(x) =
   p_B(x)$ for $x \geq j-3$, cf.\! Remark~\ref{nnn}. Since it then follows
   that $B \simeq A $ provided $B$ is Artinian, the (only) real application of
   Theorem~\ref{mainartin}(b) seems to be in the case
   $\depth_{\mathfrak{m}}B \geq 1$. It is, however, natural to use
   Theorem~\ref{mainartin}(a) also when $\depth_{\mathfrak{m}}B = 0$.
\end{remark}

\begin{remark} \label{lires}
  Suppose $\depth_{\mathfrak{m}}B \geq 1$ and that $I_{A/B}$ is
  $(j+1)$-regular, and look to
  \begin{equation}  \label{IAB}
   \small
   0 \rightarrow I_{A/B} \rightarrow B \rightarrow A \rightarrow 0 \ .
   \end{equation}
   Since $\HH_{\mathfrak{m}}^{0}(I_{A/B})=0$, we have
   $\depth_{\mathfrak{m}} I_{A/B} \geq 1$, i.e. $pd_R(I_{A/B}) \leq
   n-1$ and in fact $pd_R(I_{A/B})= n-1$ since $pd(A)= n$. A mapping
   cone construction applied to \eqref{IAB} in which we use the
   minimal resolutions of $I_{A/B}$ and $B$, leads easily to a
   $R$-free resolution of $A$. Moreover if $I_{A/B}$ admits an
   semi-linear resolution, then $(I_{A/B})_{j-1}=0$,
   and conversely provided $\reg(I_{A/B})=j+1$.  Note that $A$ becomes
   a level algebra if $I_{A/B}$ admits a linear resolution. In
   particular, the natural application of Theorem~\ref{mainartin}(b)
   is the same as for Theorem~\ref{mainzero}, cf. Remark~\ref{linres},
   i.e. the minimal resolution of $A$ should be the one of $B$ in
   addition to a semi-linear contribution coming from $I_{A/B}$ via
   the mapping cone construction.
\end{remark}
Theorem~\ref{mainartin} applies nicely to Artinian truncations and more
generally to Artinian quotients $A$ with $h$-vector $H_A=(1,h_1,
h_2,...,h_{j-1}, \alpha,0,0,..)$ where
$H_B=(1,h_1,h_2,...,h_{j-1},h_{j},h_{j+1},...)$ is the Hilbert function of $B$
and $\alpha \leq h_j$. To see it let, in a very similar way to what we did in
Proposition~\ref{sgen}, $\ \GradAlg(H)_n$ (resp. $ \GradAlg(H_B,H_A)_n$) be
the open subset of $\GradAlg(H)$ consisting of $(R/I)$ (resp. $(B =R/I
\rightarrow A)$) where the Castelnuovo-Mumford regularity satisfies $\reg(I)
\leq n$. Then we have a diagram as in \eqref{incid} where we now restrict the
natural projection morphism $ q: \GradAlg(H_B,H_A) \rightarrow \GradAlg(H_B) $
and $p$ to $\GradAlg(H_B,H_A)_n$ ;
\begin{equation} \label{introincidartin}
  \begin{array}[h]{ccc}
    \GradAlg(H_B,H_A)_n & \stackrel{q}{\longrightarrow} & \GradAlg(H_B)_n
    \subset \GradAlg(H_B) \\
    \rule[-2mm]{0pt}{6mm} \downarrow ^{p} \\
    \GradAlg(H_A)
  \end{array}
\end{equation}
 
\begin{proposition} \label{sgenartin} Let $H_B=(1,h_1, h_2, ...)$ be
  the Hilbert function of an algebra $B$ satisfying
  $\depth_{\mathfrak{m}}B \geq 1$ and let $j$, $n \leq j-2$ and
  $\alpha \leq h_j$ be integers. Let $H_A=(1,h_1, h_2,...,h_{j-1},
  \alpha,0,0,..)$ and look to the maps $p$ and $q$ in
  \eqref{incidartin}.  Then
  
  (i) \ \ $q$ is smooth and surjective with geometrically connected
  fibers, of fiber dimension $\alpha(h_j - \alpha)$, and
  
  (ii) \ \ $p$ is an isomorphism onto an open subscheme of $\ 
  \GradAlg(H_A)$. \\
  In particular the incidence correspondence \eqref{incidartin} determines a
  well-defined injective application $\pi$ from the set of irreducible
  components $W$ of $\ \GradAlg(H_B)_n$, to the set of irreducible components
  $V$ of $\ \GradAlg(H_A)$ whose general elements satify the Weak Lefschetz
  property. In this application the generically smooth components correspond.
  Indeed $V=\pi(W)$ is the closure of $p(q^{-1}(W))$, and we have
  $$
  \ \dim V = \ \dim W + \  \alpha(h_j - \alpha) \ .$$
\end{proposition}

\begin{proof}
  (i) To any point $(B')$ of $\GradAlg(H_B)_n$, let $A' :=
  \oplus_{i=0}^{j-1}B_i' \oplus V_j$ where $V_j$ is an $\alpha$-dimensional
  quotient of $B_j'$. This shows that $q$ is surjective. Moreover we get the
  smoothness of $q$ from Proposition 9(i) since $\ _0\!\HH^2(B,A,A)=0$ by
  the proof of Theorem~\ref{mainartin}.  To show that the fibers of $q$ are
  (geometrically) connected, one may look upon the fiber as the Grassmannian
  of $\alpha$-dimensional quotients of $B_j'$. Since the Grassmannian is
  irreducible, we
  conclude easily. \\
  (ii) Since the proof of the Weak Lefschetz property is standard (
  $\depth_{\mathfrak{m}}B \geq 1$), the proof is the same as for (ii) of
  Proposition~\ref{sgen}.
\end{proof}

We will call an Artinian algebra $A$ with Hilbert function $H_A$ as in
Proposition~\ref{sgenartin} with $\alpha = 0$ an {\it Artinian
  truncation in degree} $j$. Moreover, by Remark~\ref{reg+2}, we
normally need $j \geq \reg(I_B)+2$ for some $B$ to use
Proposition~\ref{sgenartin} with $ \GradAlg(H_B)_{j-2}$ non-empty.
Having several irreducible components with e.g. a non-empty
intersection in $ \GradAlg(H_B)_n$, we get exactly the same type of
irreducible components with e.g. a non-empty intersection for their
Artinian truncations in a fixed degree $j$ (for every $j \geq n+2$) in
$ \GradAlg(H_A)$ (for instance, the $B$ and the components given by
the $B_i$ of Example~\ref{twocomp}, we leave the details to the
reader). We finish this section by another example.

\begin{example} \label{nonequid} (obstructed Artinian level algebra
  with $h$-vector $(1,4,7,10,13,0,0,...)$).
  
  We have seen that $Y_1 = \Proj(B_1) \subset \pp^3$, a twisted cubic
  curve and $Y_2 = \Proj(B_2)$, a union of a plane space curve $C$ of
  degree $3$ and a point $P$ outside the plane containing $C$,
  correspond to two different irreducible components of the stratum
  $\GradAlg(H)$ of the Hilbert scheme $\Hilb^{3x+1}(\pp^3)$ where
  $H=(1,4,7,10,13,...)$. Indeed $(Y_1)$ belongs to a 12-dimensional
  irreducible component of $\GradAlg(H)$ while $(Y_2)$ belongs to
  a 15-dimensional irreducible component of $\GradAlg(H)$. Using
  Piene-Schlessinger's Theorem from \cite{PiS} to see a complete
  description of $ \Hilb^{3x+1}(\pp^3)$, we also get that the general
  element of the intersection of the two components (which is
  11-dimensional) is a curve $Y = \Proj(B)$ with an embedded point.
  The minimal resolution of $I=I_{B}$ or $I_{B_2}$ (resp. $I_{B_1}$)
  are of the form
\begin{equation} \label{resolutiB}
   0 \rightarrow R(-4) \rightarrow R(-4) \oplus R(-3)^3 \rightarrow
   R(-3) \oplus R(-2)^3  \rightarrow I \rightarrow 0 \ 
\end{equation} 
(resp. of the form \eqref{resolutiB} where both $R(-4)$ and two of
$R(-3)$ are removed). Hence the regularity of $I_B$ and $I_{B_i}$ for
$i=1$ and $2$ is at most $3$, i.e. the two components and its
intersection essentially belong to $\ \GradAlg(H)_3$. Applying
Proposition~\ref{sgenartin} for $j \geq 5$ and $n=3$ and to any
$\alpha \leq 3j+1$, we get two irreducible components $V_i$ of
$\GradAlg(H_A)$ with a well described non-empty intersection. Indeed
let $X_1= \Proj(A_1)$ and $X_2=\Proj(A_2)$ (resp. $X=\Proj(A)$) be
obtained by modding out by $h_j-\alpha$ linearly independent forms of
$(B_1)_j$ and $(B_2)_j$ (resp. $B_j$) and all forms of degree $j+1$.
It follows that $A_i$ are unobstructed as graded $R$-algebras for
$i=1$ and $2$ and that $\dim_{(A_1)} \GradAlg(H') = 12+ \alpha(h_j -
\alpha) \ $ and $\dim_{(A_2)} \GradAlg(H') = 15+ \alpha(h_j - \alpha)
\ $ where $H'=(1,4,7,...3j-5,3j-2,\alpha,0,0,...)$. Moreover $(A)$ is
a singular point of $ \GradAlg(H')$and belongs to the $11+ \alpha(h_j
- \alpha)$-dimensional intersection of the components. Finally if
$\alpha=0$ and $j=5$ it is straightforward to see that the free terms
of a minimal resolution of $A_2$ (and $A$) are
   \begin{equation*}
   \small
   0 \rightarrow R(-8)^{13} \rightarrow R(-7)^{42} \oplus R(-4) \rightarrow
   R(-6)^{45} \oplus R(-4)\oplus R(-3)^3 \rightarrow R(-5)^{16}\oplus
   R(-3) \oplus R(-2)^3 
   \end{equation*}
\end{example}

\section{Tangent and obstruction spaces of Artinian families.} 

In this section we consider graded Artinian algebras with a special
look to level quotients of $k[x,y,z]$. Note that, in most cases,
results such as Theorem~\ref{mainartin}, Proposition~\ref{sgenartin}
and Remark~\ref{summarize2} do not apply  because its assumptions limit
their applications considerably. We can, however, still
analyze $\GradAlg^{H}(R)$ at a point $(A)$ infinitesimally by means of
its tangent and obstruction spaces (and a certain obstruction
morphism, cf. \cite{L}). In the following we make these spaces more
explicit by duality, and since we show that the parameter space of
Level schemes, $\L(H)$, of \cite{GC} is essentially an open
subscheme of $\GradAlg^{H}(R)$, we can use our results to study
$\L(H)$. In particular we study in detail the level type 2
algebras which correspond to a pencil of forms by apolarity
(\cite{I04}), and we prove a conjecture of Iarrobino on the existence
of several irreducible components of $\L(H)$ when
$H=(1,3,6,10,14,10,6,2)$.

Indeed inside $\GradAlg^{H}(R)$ there is an open set consisting of graded
Artinian Gorenstein quotients $R \rightarrow A$ with Hilbert function
$H$ (which essentially is the scheme $\PGor(H)$, see \cite{K04}). An
elementary way of finding the obstruction space of $\PGor(H)$ is to
compute the kernel of the natural surjection
$$\eta_j: (S_2I_A)_j \rightarrow ({I_A}^2)_j$$
from the second
symmetric power to the second power of $I_A$ in the socle degree $j$
of $A$. Indeed, up to duality, this kernel is isomorphic to $
_0\!\HH^2(R,A,A)$, the obstruction space of $\PGor(H)$. To
generalize this result to any Artinian $A$, we remark that $\ker
\eta_j$ is isomorphic to the cokernel of the natural morphism $
(\Lambda_2 I_A)_j \rightarrow \Tor^R_2(A,A)_j$ (at least if $char(k)
\neq 2$). This formulation allows a generalization to any Artinian
$A$. Indeed there is a special product, given as an antisymmetrization
map (\cite{AND}, Prop.\! 24.1),
\begin{equation} \label{TorK}
 \Tor^R_1(A,A) \otimes_A \Tor^R_1(A,K_A) \rightarrow \Tor^R_2(A,K_A)
 \ 
\end{equation}
with cokernel $\HH_2(R,A,K_A)$. Up to duality we will show that the zero
degree piece of this cokernel is the obstruction space of $\GradAlg^{H}(R)$ at
$(A)$. To prove it we need a variation of the following spectral sequence
  \begin{equation*} \label{spectalg}
   \Ext_A^p(\HH_q(R,A,A),M)   \implies \ \HH^{p+q}(R,A,M) \ 
 \end{equation*}
 (cf.\! \cite{AND}, Prop. 16.1). Keeping also the spectral sequence
 (\cite{HER}, Satz\! 1.2)
 \begin{equation} \label{spectext}
 \Ext_C^p(\Tor^A_q(M,K_C),K_C) \implies \ \Ext_A^{p+q}(M,C)
 \end{equation}
 in mind ($C$ a CM quotient of $A$ with canonical module $K_C$), the following
 result is not surprising

\begin{proposition} \label{spectprop}
If $B \rightarrow A \rightarrow C$ are quotients of $R$ of arbitrary
  dimension and if $C$ is CM with canonical module $K_C$, then there
  is a spectral sequence converging to $  \HH^{*}(B,A,C)$ ;
 \begin{equation} \label{spectalgK}
   'E_2^{p,q}:=   \Ext_C^q(\HH_p(B,A,K_C),K_C)  \implies \
   \HH^{p+q}(B,A,C) \ .
 \end{equation}
 In particular if $C$ is a graded Artinian algebra, then there is a
 graded preserving isomorphism $$
 \Hom_C(\HH_q(B,A,K_C),K_C) \simeq \ 
 \HH^{q}(B,A,C) \ .$$
\end{proposition}

\begin{proof} One knows that $ \Hom_A(M,C) \simeq   \Hom_C(M \otimes_A
  K_C,K_C)$, $M$ an $A$-module. Using this we can prove our
  proposition in the usual way, i.e. by considering the double complex
  $$K_{*,*}= \Hom_C(\Diff(B,A_*,A) \otimes_A K_C, I_*)$$
  where $0
  \rightarrow K_C \rightarrow I_*$ is an injective resolution of the
  $C$-module $K_C$ and $\Diff(B,A_*,A):= \Omega_{A_*/B} \otimes_{A_*} A
  $ is the complex of K\"ahler differentials based on a simplicial
  resolution, $A_*$, of the $B$-algebra $A$ (as in \cite{AND}, Prop.\!
  17.1, so each $A_n$ is a
  polynomial ring over $B$). 
  If we in $''E_2^{p,q}$ first take homology of $K_{*,*}$ with respect to the
  second variable (i.e. $I_*$), we get $''E_2^{p,0}= \HH^{p}(B,A,C)$ and
  $''E_2^{p,q}=0$ for $q >0$ because $ \Ext_C^q(K_C,K_C)=0$ for $q>0$ by
  Cohen-Macaulayness and the fact that each $\Diff(B,A_n,A) \otimes_A C$ is
  $C$-free. Calculating $'E_2^{p,q}$ by reversing the order, i.e. by first
  taking homology with respect to the first variable, we get
  \eqref{spectalgK}.  Finally since $K_C$ is an injective $C$-module in the
  Artinian case, we are done.
\end{proof}

\begin{theorem}\label{mainartingrad}
  Let $R \rightarrow A = R/I_A$ be a graded Artinian quotient with
  Hilbert function $H$. Then the dual of $ (I_A \otimes_R K_A)_0$ is
  the tangent space of $ \GradAlg^{H}(R)$ at $(A)$, and the dual of $
  _0\!\HH_2(R,A,K_A)$ contains the obstructions of deforming $A$ as a
  graded $R$-algebra. Moreover $$\dim\ (I_A \otimes_R K_A)_0 - \ 
  _0\!\hho_2(R,A,K_A) \leq \dim_{(A)}\! \GradAlg(H) \leq \dim\ (I_A
  \otimes_R K_A)_0.$$
  In particular $ \GradAlg^{H}(R)$ is smooth at
  $(A)$ 
  provided the natural map $$I_A \otimes_R I_A \otimes_R K_A
  \rightarrow \Tor^R_1(I_A,K_A)$$
  of \eqref{TorK}, i.e. the map
  $\zeta$ concretely described in \eqref{TorIK} below, is surjective in
  degree zero. 
\end{theorem}

\begin{proof}
  Since it is known that the tangent (resp.\ ``obstruction'') space of
  $ \GradAlg(H_A)$ at $(A)$ is $_0\!\HH^1 (R,A,A)= \ 
  _0\!\Hom_A(I_{A}/I_A^2,A) \simeq \ _0\!\Hom_R(I_{A},A) \ $ (resp.\ 
  $_0\!\HH^2 (R,A,A)$) by \cite{K79}, Thm.\! 1.5, we get the
  description in Theorem~\ref{mainartingrad} of these spaces by
  Proposition~\ref{spectprop}. Then the left inequality of the
  dimension formula follows rather easily from \cite{L}, Thm.\! 4.2.4
  while the right inequality is trivial. Hence we get all conclusions
  of the theorem once we have shown the surjectivity in \eqref{TorK}
  and the surjectivity of $\zeta$ in \eqref{TorIK} are equivalent.
  Indeed $ \Tor^R_2(A,K_A) \simeq \Tor^R_1(I_A,K_A)$ and $
  \Tor^R_1(A,K_A) \simeq I_A \otimes_R K_A$ and the map of
  \eqref{TorK} is just the natural map $\zeta: I_A \otimes_R I_A
  \otimes_R K_A \rightarrow \Tor^R_1(I_A,K_A)$ uniquely described in
  the following way. Let $0 \rightarrow N \rightarrow F \rightarrow
  K_A \rightarrow 0$ be a short exact sequence where $F$ is $A$-free.
  Applying $I_A \otimes_R (-)$ onto this sequence we get an injection
  $ \Tor^R_1(I_A,K_A) \hookrightarrow I_A \otimes N$ which together
  with the surjection $ F \twoheadrightarrow K_A $ lead to the
  composition
  \begin{equation} \label{TorIK} 
I_A \otimes_R I_A \otimes_R F \twoheadrightarrow 
  I_A \otimes_R I_A \otimes_R K_A \stackrel{\zeta}{\longrightarrow}
  \Tor^R_1(I_A,K_A) \hookrightarrow I_A \otimes_R N 
   \end{equation}
   given by $ x \otimes y \otimes \omega \rightarrow x \otimes (y \omega) - y
   \otimes (x \omega) \in I_A \otimes_R N$ (cf. \cite{AND2}, Prop.\! 9, p.\!
   204 for details).
\end{proof}

\begin{remark} \label{remHnede}
  Let $M$ be an graded $A$-module and let $0 \rightarrow N \rightarrow
  F \rightarrow M \rightarrow 0$ be a graded exact sequence where $F$
  is $A$-free.  Arguing as in the proof above, we see that $
  _v\!\HH_2(R,A,M)$ is the homology in degree $v$ of \eqref{vanHnede}
  below. Hence it vanishes if and only if the sequence
  \begin{equation}  \label{vanHnede} 
   I_A \otimes_R I_A \otimes_R F \stackrel{\lambda}{\longrightarrow} I_A
   \otimes_R N \rightarrow I_A \otimes_R F \ , 
   \end{equation}
   where $ \lambda(x \otimes y \otimes \omega)= x \otimes (y \omega) -
   y \otimes (x \omega)$, is exact in degree $v$ (\cite{AND2}, Prop.\!
   9, p.\!  204).
\end{remark}

\begin{remark} \label{remspektr} Let $A = R/I_{A}$ be a graded
  Artinian algebra and let $M$ be a finitely generated $R$-module. Using
  \eqref{spectext} we get $$
  \Hom_A(\Tor^R_q(M,K_A),K_A) \simeq \ 
  \Ext_R^{q}(M,A) \ .$$
  Thus $ (I_A \otimes_R K_A)_v$ (resp. $_v\!
  \Tor^R_1(I_A,K_A)$) is dual to $_{-v}\!  \Hom_R(I_A,A)$ (resp. $_{-v}\!
  \Ext_R^{1}(I_A,A)$) and the dual of the degree $v$ part of
  \eqref{TorK} augmented by $ \ _v\!\HH_2(R,A,K_A)$ 
 yields an exact sequence  
 $$
 \ _{-v}\!\HH^2(R,A,A) \hookrightarrow \ _{-v}\!\Ext_R^1(I_A,A)
 \rightarrow \ _{-v}\!\Hom_R(I_A \otimes_R I_A ,A)$$
 where the left
 injective map must be the
right inclusion of \eqref{cohom} in degree $-v$.
\end{remark} 

In the codimension 3 case it turns out that $\ _{-v}\!\Ext_R^1(I_A,A)$ is also
dual to $\ _{v-3}\! \Hom_R(I_A,A)$:
\begin{proposition} \label{euler} Let $R \rightarrow A = R/I_A$,
  $R=k[x,y,z]$ be a graded Artinian quotient with Hilbert function $H$
  and minimal resolution
  \begin{equation} \label{resA} 0 \rightarrow \oplus_{i=1}^{r_3} R(-n_{3,i})
    \rightarrow \oplus_{i=1}^{r_2} R(-n_{2,i}) \rightarrow \oplus_{i=1}^{r_1}
    R(-n_{1,i})\rightarrow R \rightarrow A \rightarrow 0 \ .
\end{equation}
Then $_v\! \Ext_R^{i}(I_A,A) \simeq \ _{v}\! \Tor^R_{1-i}(I_A,K_A(3))$ for $0
\leq i \leq 1$ and $N_A:= \Hom_R(I_A,A)$ satisfies $$\dim (N_A)_v - \
_v\!\ext_R^1(I_A,A)= \sum_{j=1}^3
\sum_{i=1}^{r_j}(-1)^{j-1}H(n_{j,i}+v)-H(-v-3) \ .$$ In particular $_v\!
\Ext_R^{1}(I_A,A)$ is dual to $(N_A)_{-v-3}$ for every $v$. Moreover if
$(N_A)_{-3}=0$, then $_{0}\!\HH^2(R,A,A)=0$ and $\GradAlg(H)$ is smooth at
$(A)$ of dimension $ \sum_{j=1}^3 \sum_{i=1}^{r_j}(-1)^{j-1}H(n_{j,i})$.
\end{proposition}

\begin{proof} Dualizing  \eqref{resA} we get an $R$-free resolution of
  $K_A(3)$. Moreover we essentially get the complex
 \begin{equation*} \label{resAdual} 0 \rightarrow \oplus_{i=1}^{r_1} A(n_{1,i})
    \rightarrow \oplus_{i=1}^{r_2} A(n_{2,i}) \rightarrow \oplus_{i=1}^{r_1}
    A(n_{3,i})\rightarrow 0 \ 
\end{equation*}
by either tensoring the resolution of $K_A(3)$ by $A$, or by applying $_0\!
\Hom_R(-,A)$ onto the minimal resolution of $I_A$ deduced from \eqref{resA}.
In particular $_v\! \Ext_R^{i}(I_A,A) \simeq\ _{v}\! \Tor^R_{2-i}(A,K_A(3))
\simeq \ _{v}\! \Tor^R_{1-i}(I_A,K_A(3))$ for $0 \leq i \leq 1$ and $_v\!
\Ext_R^{2}(I_A,A) \simeq\ K_A(3)_{v}$ whose dimensions satisfy the doubly
summation formula. We conclude by Remark~\ref{remspektr} and
Theorem~\ref{mainartingrad}.
\end{proof}

\begin{corollary}  With $A$ as in Proposition~\ref{euler} we have
  the equality to the right
  $$
  \rho(H) := \sum_{j=1}^3 \sum_{i=1}^{r_j}(-1)^{j-1}H(n_{j,i})= 1 -
  \sum_{j=1}^3 \sum_{i=1}^{r_j}(-1)^{j-1}H(n_{j,i}-3) \ . $$
  Moreover the sums
  above, which we call $\rho(H)$, depends only upon the Hilbert function H and
  not upon the graded Betti numbers. We have $$\dim_{(A)}
  \GradAlg^H(k[x,y,z]) \geq \rho(H)\ .$$
  In particular $\rho(H)$ is a lower
  bound for the dimension of any irreducible component of $
  \GradAlg^H(k[x,y,z])$. 
\end{corollary} 

\begin{proof} The duality of the proposition shows that $\dim (N_A)_v - \
  _v\!\ext_R^1(I_A,A)= \ _{-v-3}\!\ext_R^1(I_A,A)- \dim (N_A)_{-v-3}$. Putting
  $v=0$ we get the equality of the two expressions of $\rho(H)$ of the
  corollary from Proposition~\ref{euler}. Moreover, by \cite{L}, Thm.\! 4.2.4,
  the number $\dim (N_A)_0 - \dim \ _{0}\!\HH^2(R,A,A)$ is a lower bound of
  $\dim_{(A)} \GradAlg^H(k[x,y,z])$. Using Remark~\ref{remspektr} we get
  $\rho(H)$ to be a possibly smaller lower bound.  Finally the sum which
  defines $ \rho(H)$ depends only upon the Hilbert function because the
  contribution from all ghost terms sums to zero!
\end{proof}
\begin{example} \label{exlink} To illustrate Proposition~\ref{euler},
  we consider $H=(1,3,6,6,3,1)$ and the two different irreducible
  components (now of $ \GradAlg^H(k[x,y,z])$) of
  Example~\ref{ekslicci} whose general elements are licci. Looking to
  the minimal resolutions of $A_i$ of Example~\ref{ekslicci}, we get
  $$\ \dim (N_{A_i})_0 - \ _0\!\ext_R^1(I_{A_i},A_i)= 4H(3)-4H(5)=20 \ .$$
  By
  Remark~\ref{remspektr}, $\ _0\!\Ext_R^1(I_{A_1},A_1)=0$ since $\ _{0}\!
  \Tor^R_{1}(I_{A_1},K_{A_1}) \simeq (I_{A_1} \otimes I_{A_1})_5 =0$. (Indeed
  this $_0\!\Ext^1$ -group always vanishes in the compressed Gorenstein case.)
  Thus the ``Gorenstein'' component has dimension $20$ (also well known by
  \cite{I84}), while Remark~\ref{summarize2}(iii) or Theorem~\ref{mainlink}
  applied to the successive linkages $(1,2,3)$, $(2,2,4)$, $(2,3,4)$,
  $(3,4,4)$, $(3,4,5)$, $(3,3,5)$ obtained from a CI of type $(1,1,1)$, shows
  that the other component is generically smooth of dimension $21$. Thus $
  _0\!\ext_R^1(I_{A_2},A_2)=1$.
\end{example}

Recall that if $A$ itself admits a semi-linear $R$-free resolution
(except possibly at the minimal generators of $I_A$), then
\begin{equation} \label{vanishH2}
_0\!\HH^2(R,A,A)= 0 \  
\end{equation}
by Remark~\ref{H2=0} and Remark~\ref{lires}. This vanishing also follows
from Theorem~\ref{mainartingrad}. Moreover using
Theorem~\ref{mainartingrad}, we can prove a ``dual'' result. Indeed
suppose $I_A$ admits a semi-linear resolution except possibly at
the {\it left} end of the resolution, i.e. suppose $I_A$ has minimal
generator only in degree $j$ and $j+1$ and that the resolution
continues by
\begin{equation} \label{linearI}
  0 \rightarrow G \oplus R(-j-n+1)^{\alpha_n} \rightarrow
  R(-j-n+1)^{\beta_{n-1}} \oplus R(-j-n+2)^{\alpha_{n-1}} \rightarrow
  . .  .  \rightarrow F_1 \rightarrow I_A 
\end{equation}
where $G$ is any $R$-free module. Here $F_1= R(-j-1)^{\beta_1} \oplus
R(-j)^{\alpha_1}$, $R$ is the polynomial ring $k[x_1,...,x_n]$ and $n\geq 3$.
Then \eqref{vanishH2} holds. Using \eqref{spectext} we can even replace $F_1$
by
\begin{equation} \label{linearcorr} F_1= R(-j-1)^{\beta_1} \oplus
  R(-j)^{\alpha_1} \oplus ( \oplus_{i=1}^m R(-a_i)) \ \ , \ \ a_i < j \ \ {\rm
    for \ all} \ \ i
\end{equation}
where the set of generators $ \{f_1,...,f_m \} $ which correspond to $
\{a_1,...,a_m \} $ form a regular sequence, and still get \eqref{vanishH2},
i.e.

\begin{proposition} \label{corart} Let $A = R/I_A$ be a graded Artinian
  quotient with Hilbert function $H$, whose minimal resolution is given by
  \eqref{linearI} and \eqref{linearcorr} where the generators $ \{f_1,...,f_m
  \} $ of $I_A$ which correspond to $ \{a_1,...,a_m \} $ form a regular
  sequence. Let $B = R/(f_1,...,f_m )$ (and $B=R$ if $m=0$). Then $\
  _0\!\HH^2(R,A,A)= 0$ and $ \GradAlg(H)$ is smooth at $(A)$. Moreover $$
  \dim_{(A)}\! \GradAlg(H) = \ _{-n}\! \hom_R(G,B) - _{-n}\! \hom_R(G,A)+
  \sum_{i=1}^m H(a_i) . $$ 
 \end{proposition}

\begin{proof}
  By the long exact sequence of algebra cohomology (Remark~\ref{H2=0}) and
  \eqref{cohom} we get $\ _0\!\HH^2(R,A,A)= 0$ provided we can show $
  \Ext_B^{1}(I_{A/B},A)=0.$ Continuing the long exact sequence of
  Remark~\ref{H2=0} to the left we see that $ \Ext_B^{1}(I_{A/B},A)=0$ also
  leads to
  \begin{equation} \label{dimen} \dim (N_A)_0= \! \hom_B(I_{A/B},A) +
    \hom_B(I_B/I_B^2,A) \ .
\end{equation} 
To show $ \Ext_B^{1}(I_{A/B},A)=0$, we improve a little bit upon
Theorem~\ref{mainartingrad} by using \eqref{spectext}. Indeed we have $
\Hom_A(\Tor^B_q(I_{A/B},K_A),K_A) \simeq \ \Ext_B^{q}(I_{A/B},A)$. Hence it
suffices to show $ _0\! \Tor^R_1(I_{A/B},K_A)=0$. Now look to the exact
sequence
\begin{equation*} \label{linearK} \rightarrow R(j+n-1)^{\beta_{n-1}} \oplus
  R(j+n-2)^{\alpha_{n-1}} \rightarrow G^* \oplus R(j+n-1)^{\alpha_n}
  \rightarrow K_A(n) \rightarrow 0
\end{equation*} 
which we tensorize with $I_{A/B}(-n)$. By the definition of
$\Tor^R_1(I_{A/B},K_A)$, it suffices to show $ (I_{A/B}(-n)(j+n-1))_0
= 0$ and $(I_{A/B}(-n)(j+n-2))_0=0$. This is true since
$(I_{A/B})_{j-1}=0$ by assumption. Moreover the argument also shows
$(I_{A/B} \otimes K_A )_0 \simeq \dim (G^*(-n) \otimes I_{A/B})_0$.
Hence we get $$_0\! \hom_B(I_{A/B},A) = \dim \Tor^B_0(I_{A/B},K_A)= \ 
_{-n}\! \hom_R(G,B) -\ _{-n}\! \hom_R(G,A) \ ,$$
and we conclude by
\eqref{dimen} and the fact that $I_B/(I_B)^2 \simeq \oplus_{i=1}^m
B(-a_i)$.
\end{proof} 

Proposition~\ref{corart} with $B=R$ applies nicely to compressed Artinian
algebras. Indeed the number $ \ _{-n}\! \hom_R(G,B) - _{-n}\! \hom_R(G,A)$
coincides with the dimension of the corresponding component given in Thm.\!
IIB of \cite{I84}. If $B \neq R$ Proposition~\ref{corart} also applies to
non-compressed algebras:

\begin{example}
  As a special case of Proposition~\ref{corart} we look to Artinian level
  quotients with ``almost semi-linear'' resolution. All level algebras below
  may be constructed as $A=R/ann(F_1,F_2)$ where $F_1$ and $F_2$ are forms of
  the degree $7$ in the dual polynomial algebra of $R$ (cf. later discussion).
  Indeed we easily construct in this way algebras $A_i$ with Hilbert functions
  $H_{A_1}=(1,3,6,10,15,12,6,2)$, $H_{A_2}=(1,3,6,10,14,12,6,2)$,
  $H_{A_3}=(1,3,6,10,13,12,6,2)$ 
  and corresponding minimal resolutions
\begin{equation*} \label{resolLevel1}
   0 \rightarrow R(-10)^2 \rightarrow R(-7)^5 \oplus R(-6)^5 \rightarrow
   R(-5)^9   \rightarrow I_{A_1} \rightarrow 0 \ 
\end{equation*} 
\begin{equation*} \label{resolLevel2}
   0 \rightarrow R(-10)^2 \rightarrow R(-7)^6 \oplus R(-6)^4 \rightarrow
   R(-6)^2 \oplus R(-5)^6 \oplus R(-4)  \rightarrow I_{A_2} \rightarrow 0 \ 
\end{equation*} 
\begin{equation*} \label{resolLevel3}
   0 \rightarrow R(-10)^2 \rightarrow R(-7)^7 \oplus R(-6)^3 \rightarrow
   R(-6)^4 \oplus R(-5)^3 \oplus R(-4)^2  \rightarrow I_{A_3} \rightarrow 0 \ .
\end{equation*} 
Only $A_1$ is compressed, but since one may show that the minimal generators
of $I_{A_i}$ of degree $4$ (which we use to define $B_i$) form a regular
sequence,
Proposition~\ref{corart} applies (we have used Macaulay 2 to check it and to
find the minimal resolutions). Hence the algebras $A_i$ are unobstructed and
since $ \ _{-n}\! \hom_R(G,M) = 2 \cdot \dim M_7$ and $\dim(B_i)_7 = (i-1)\dim
R_3 $ for $i=2$ and 3 (and $B_1=R$), we get the number $$\ \dim_{(A_i)}\!
\GradAlg(H_{A_i})= 2 \cdot \dim (B_i)_7 - 2\cdot \dim H_{A_i}(7) +
(i-1)H_{A_i}(4)$$
to be $68, 62, 54$ for $i=1,2,3$ respectively.
\end{example}
To this end we consider level algebras of CM-type $t$. Let $\LevAlg(H)$ be the
open set of $\GradAlg(H)$ (and hence open as a subscheme with its induced
scheme structure) consisting of graded level quotients with Hilbert function
$H$. Since we work with Artinian algebras there is another known scheme,
$\L(H)$, parametrizing graded level quotients with suitable Hilbert function
$H$, namely the determinantal loci in the Grassmannian $G(t,j)$ of
$t$-dimensional vector spaces of forms of degree $j$, cut out by requiring
their ``catalecticant matrices'' to have ranks given by the Hilbert function
(see \cite{GC}, and \cite{IK}, Sect.\! 1.1 for the Gorenstein case). Then the
underlying sets of closed points of $\L(H)$ and $\LevAlg(H)$ are the same by
apolarity (the Macaulay correspondence), and their tangent spaces are
isomorphic (\cite{GC}, Thm.\! 2.1 for $\L(H)$, and \cite{K79}, Thm.\! 1.5 for
$\GradAlg(H)$). Since one may by the proof below see that $\LevAlg(H)$ and
$\L(H)$ are in fact isomorphic as topological spaces (expected since they have
the Zariski topologies and the bijection between them is natural), we have

\begin{theorem}\label{mainlev}
  Let $R \rightarrow A$ be a graded Artinian level quotient with Hilbert
  function $H$. Then $\dim _{(A)} \GradAlg^{H}(R) = \dim _{(A)} \L(H)$. Hence
  $\L(H)$ is smooth at $(A)$ if and only if $\GradAlg^{H}(R)$ is smooth at
  $(A)$. In particular $\L(H)$ is smooth at $(A)$ provided $ _0\!\HH^2(R,A,A)=
  0$, i.e. provided the map of \eqref{TorK} is surjective in degree zero, or
  equivalently, the displayed sequence of Remark~\ref{remHnede} with $M=K_A$
  is exact in degree zero.
\end{theorem}

\begin{proof} Let $V \subset \L(H)$ be a closed irreducible subset, 
  and let $V$ have the reduced scheme structure. By the definition of $\L(H)$,
  the restriction of the ``universal'' bundle of the Grassmannian $G(t,j)$ to
  $V$ defines via apolarity a family of graded Artinian level quotients,
  $A_V$, over $V$ with constant Hilbert function $H$. Since $V$ is integral,
  it follows that the family (i.e. the morphism $\Spec(A_V) \rightarrow V$) is
  flat (\cite{M}, Lect.\! 6). Hence we have a morphism $\pi: V \rightarrow
  \LevAlg(H) $ by the universal property of $\GradAlg(H)$.  $\pi(V)$ is
  irreducible and closed in $\LevAlg(H)$ (it is closed because an ``inverse''
  $(\LevAlg(H))_{red} \rightarrow \L(H)$ on closed points exists, by
  \cite{IK}, p. 249). So chains of closed irreducible subsets in $\L(H)$ and
  $\LevAlg(H)$ correspond, and the spaces have the same dimension. Since their
  tangent spaces are isomorphic, it follows that $\GradAlg^{H}(R)$ is smooth
  at $(A)$ iff $\L(H)$ is smooth at $(A)$. Then we conclude by
  Theorem~\ref{mainartingrad} since the surjectivity of \eqref{TorK} in degree
  zero is equivalent to the exactness of the corresponding sequence in
  Remark~\ref{remHnede}.
\end{proof}

As an application we consider certain type 2 level algebras studied by
Iarrobino in \cite{I04}, i.e. level algebras given by $A=R/ann(F_1,F_2)$ where
$F_1$ and $F_2$ are forms of the same degree $j$ in the dual polynomial
algebra of $R$, upon which $R$ acts by differentiation. Let
$A_i:=R/ann(F_i)$. Since we have $I_A = I_{A_1} \cap I_{A_2}$, we get an exact
sequence
$$0 \rightarrow A \rightarrow A_1 \oplus A_2 \rightarrow R/(ann(F_1)+ann(F_2))
\rightarrow 0.$$
In an extended draft of \cite{I04} the author determines the
tangent space of $ \LevAlg(H)$ at such an $(A)$ and he gives it a manageable
form in the case $\{F_1,F_2\}$ is {\it complementary}, i.e. provided
\begin{equation}  \label{complement}
H_A(i)= \min \{\dim R_i, H_{A_1}(i)+H_{A_2}(i) \} \ {\rm \ \ for \ any} \ i \ .
\end{equation} 
where $H_A=H$. Our Theorem~\ref{mainartingrad} gives us, not only a
tangent space which coincides with his, but it provides us also with
the following manageable form of the obstruction space.

\begin{proposition}\label{mainlev2}
  Let $\{F_1,F_2\}$ be complementary forms of degree $j$, and let
  $A = R/I_A$ be the Artinian level quotient with Hilbert
  function $H$ given by $I_A=ann(F_1,F_2)$. Let $I_{A_i}=ann(F_i)$.
  Then $(I_A/I_A \cdot I_{A_1})_j\ \oplus (I_A/I_A \cdot I_{A_2})_j$
  is the dual of the tangent space of $\ \GradAlg^{H}(R)$ at $(A)$, and
  $ _j\! \HH_2(R,A,A_1) \oplus \!  _j\!\HH_2(R,A,A_2) $ is the dual of a
  space containing all obstructions of deforming $A$ as a graded
  $R$-algebra.  In particular if the sequences
\begin{equation*}  
   I_A \otimes_R I_A  \stackrel{\lambda}{\longrightarrow} I_A
   \otimes_R I_{A_i} \twoheadrightarrow I_A \cdot I_{A_i}   
   \end{equation*}
   where $ \lambda(x \otimes y)= x \otimes y - y \otimes x$, are exact in
   degree $j$ for $i=1,2$, then $\GradAlg(H)$ (and $\L(H)$) is smooth at $(A)$
   and we have
  $$
  \ \dim_{(A)} \GradAlg^{H}(R)= \ \sum_{i=1}^2 \dim (I_A/I_A \cdot
  I_{A_i})_j 
  $$
\end{proposition}

\begin{remark} \label{remmainlev}
  The map $I_A \otimes_R I_A \stackrel{\lambda'}{\longrightarrow} I_A
  \otimes_R I_{A}$, defined by $ \lambda'(x \otimes y)= x \otimes y -
  y \otimes x$, obviously commutes with $\lambda$ above. Since
  $\lambda'$ factors via the natural surjection $I_A \otimes_R I_A
  \twoheadrightarrow \Lambda^2I_A$ (in $char(k) \neq 2$), then
  $\lambda$ also does. In $char(k) \neq 2$ the exactness of the two
  sequences of Proposition~\ref{mainlev2} is therefore equivalent to
  the exactness of
\begin{equation} \label{exter} 
 \Lambda^2I_A \rightarrow I_A \otimes_R I_{A_i}
  \twoheadrightarrow I_A \cdot I_{A_i} \ , 
\end{equation} 
$i=1,2$, in degree $j$. Indeed, by Remark~\ref{remHnede}, $ _j\!
\HH_2(R,A,A_i)$ is the homology of \eqref{exter} in degree $j$. In particular
if $(I_A \otimes_R I_A)_j \simeq (S_2I_A)_j$ (e.g. $(I_A \otimes_R I_A)_j
\simeq ({I_A}^2)_j$), then the exactness of the sequences of
Proposition~\ref{mainlev2} is equivalent to $(I_A \otimes_R I_{A_i})_j \simeq
(I_A \cdot I_{A_i})_j$.
\end{remark} 

\begin{proof}
  Let $s(I_A)$ be the minimal degree of a minimal generator of $I_A$ and let
  $\overline{A} =R/(ann(F_1)+ann(F_2)) \ $. Since $\{F_1,F_2\}$ is
  complementary, we get $(\overline{A})_v = 0$, i.e. $A_v \simeq (A_1)_v
  \oplus (A_2)_v$ for $v \geq s(I_A)$. It follows that $$ (K_{A_1})_v \oplus
  (K_{A_2})_v \simeq (K_{A})_v$$ for $v \leq -s(I_A)$. Defining $\overline{K}$
  by the long exact sequence 
\begin{equation}  \label{barK}
  0 \rightarrow \overline{K} \rightarrow K_{A_1} \oplus
  K_{A_2} \rightarrow K_{A} \rightarrow 0 \ ,
\end{equation}
we get  $(\overline{K})_v = 0$  for $v \leq -s(I_A)$. By 
considering a minimal $R$-free resolution of $I_A$, it follows that 
\begin{equation}  \label{torbarK}
  _0\! \Tor_i^R(I_A, \overline{K}) = 0 \ \ {\rm for} \ \ i \leq 0 \ .
   \end{equation}
Now applying $I_A \otimes (-)$ onto \eqref{barK}, or more precisely using the
corresponding long exact sequence of algebra homology, we get $$
_0\! \HH_i(R,A,K_{A_1} \oplus K_{A_2}) \simeq \  _0\! \HH_i(R,A, K_{A})
$$
for $i=1$ and $2$ because $ _0\! \Tor_1^R(I_A, \overline{K}) \simeq\ _0\!
\Tor_2^R(A, \overline{K}) \twoheadrightarrow \ _0\! \HH_2(R,A, \overline{K})$
is surjective (cf. \eqref{TorK}) and $I_A \otimes_R \overline{K} \simeq
\HH_1(R,A, \overline{K})$, i.e. $ _0\! \HH_i(R,A, \overline{K})$ vanishes for
$i=1$ and $2$ by \eqref{torbarK}. Then we conclude easily by $A_i \simeq
K_{A_i}(-j) $, Theorem~\ref{mainartingrad} and Remark~\ref{remHnede}. Indeed
we have $$(I_A \otimes_R K_{A_i})_0 \simeq \ (I_A \otimes_R A_i(j))_0 \simeq
(I_A \otimes_R R/I_{A_i})_j \simeq (I_A/I_A \cdot I_{A_i})_j $$ and we get
$(I_A \otimes_R K_{A})_0 \simeq \ (I_A \otimes_R K_{A_1})_0 \oplus (I_A
\otimes_R K_{A_2})_0 $ as well as
$$
_0\! \HH_2(R,A,K_A) \simeq \ _j\! \HH_2(R,A,A_1) \oplus \!
_j\!\HH_2(R,A,A_2) \ . $$
By the assumption of the exactness of the sequences
and by Remark~\ref{remHnede} (letting $F=R$ and $N=I_{A_i}$), we get the
vanishing of $ _0\! \HH_2(R,A,K_A) \ $ and we are done.
\end{proof}

\begin{remark} \label{remcomp} As Iarrobino points out in the draft
  of \cite{I04}, Theorem $4.8A$ of \cite{I84} implies that if
  $F_1$ is any form of degree $j$ and $F_2$ is a sum of length $s$ of
  linear forms to the $j$-th power (i.e. the Hilbert function of
  $A_2=R/ann(F_2)$ equals $H(s,j,n)$ of Remark~\ref{summarize2}(v)),
  then $\{F_1,F_2\}$ is complementary provided we choose the linear
  forms of $F_2$ general enough. It follows that $H_A$ is given by
  \eqref{complement}.
\end{remark}

Firstly we give an easy example which may also be treated by
Proposition~\ref{corart}. 

\begin{example} \label{excomp} 
  
  (a) Let $H=(1,3,6,10,6,2)$. By Remark~\ref{remcomp} there are forms
  $F_1$ and $F_2$ where each $F_i$ is a sum of length $5$ of linear
  forms to the $5$-th power (i.e. $H_{A_i}=(1,3,5,5,3,1)$) and such
  that the Hilbert function of $A=R/ann(F_1,F_2)$ is $H$. Then $A$ is
  compressed. From the Hilbert functions we see that $s(I_A)=4$ while
  $s(I_{A_i})=2$. Moreover, the socle degree of $A$ and $A_i$ are $5$,
  and we get $(I_A \otimes I_{A_i})_5 = 0$ for $i=1$ and $2$. By
  Proposition~\ref{mainlev2} it follows that $ \GradAlg(H)$ is
  unobstructed at $(A)$ and we have
  $$
  \ \dim_{(A)} \GradAlg(H)= \ 2 \cdot \dim (I_A)_5 = 38 \ .$$

  (b) Let $H=(1,3,6,9,6,2)$, let $F_1$ be as in (a), while we now let $F_2$
  be a sum of length $4$ of general linear forms to the $5$-th power. Hence
  $H_{A_2}=(1,3,4,4,3,1)$ and $H = H_A$ where $A=R/ann(F_1,F_2)$ by
  Remark~\ref{remcomp}. From the Hilbert functions we see that $s(I_A)=3$ and
  $s(I_{A_i})=2$. Since we easily see that $I_A \otimes I_{A_i} \simeq I_A
  \cdot I_{A_i}$ is an isomorphism in degree $5$ for $i=1$ and $2$, we get by
  Proposition~\ref{mainlev2} that $ \GradAlg(H)$ is unobstructed at $(A)$
  and that 
  $$
  \ \dim_{(A)} \GradAlg(H)= \ 2 \cdot \dim (I_A)_5 - dim (I_A \cdot
  I_{A_1})_5 - \dim (I_A \cdot I_{A_2})_5 = 35 \ .$$
\end{example}

Loosely speaking it is, for $i=1,2$, the relations of $I_A \cdot I_{A_i}$ in
degree $j$, modulo those coming from the relations of $I_A \otimes I_{A_i}$
and the generators of $\wedge^2 I_A$, which contribute to $ _0\!
\HH^2(R,A,A)$. Of course the vanishing of $ _0\! \HH^2(R,A,A)$ as well as the
dimension of $ \ \GradAlg^{H}(R)$ is usually straightforward to get from
Proposition~\ref{mainlev2} provided $s(I_A)+s(I_{A_i}) \geq j$ for $i=1,2$, as
in Example~\ref{excomp}.
%

We finish this paper by proving a conjecture of Iarrobino, appearing in the
draft of \cite{I04}, namely that $\L(H)$ with $H=(1,3,6,10,14,10,6,2)$
contains at least two irreducible components, where one of the components
contains Artinian level type 2 algebras given by 2 forms of Hilbert function
$H_1=(1,3,6,9,9,6,3,1)$ and $H_2=(1,3,4,5,5,4,3,1)$, as in
Remark~\ref{remcomp}, and the other contains level type 2 algebras constructed
via 2 forms with Hilbert function $H_3=(1,3,6,10,10,6,3,1)$ and
$H_4=(1,3,4,4,4,4,3,1)$. As pointed out in the Introduction, even though this
conjecture was open until now, Iarrobino and Boij have in a joint work already
constructed other examples of reducible $\L(H)$ whose general elements are
level quotients of type 2, one with $H=(1,3,6,10,14,18,20,20,12,6,2)$, and
moreover got a doubly infinite series of such components \cite{BI}.

\begin{example} \label{lev2comp} Let $H=(1,3,6,10,14,10,6,2)$. We claim that
  there are at least two components $V_1$ and $V_2$ of $\L(H)$ whose general
  elements are Artinian level type 2 algebras, that $\dim V_1 = 46$ and $\dim
  V_2 = 47$ and that both components are generically smooth.
   
  To get the component $V_1$ of dimension $46$, take $F_1$ to be a sum of
  length $4$ of general linear forms to the $7$-th power and take
  $F_2$ to be a general polynomial of degree 7. If $A_i=R/ann(F_i)$
  and $A=R/ann(F_1,F_2)$, then $H_{A_2}=(1,3,6,10,10,6,3,1)$,
  $H_{A_1}=(1,3,4,4,4,4,3,1)$ and $H_A=H$. It suffices to show that
  $A$ is unobstructed and that $ \ \dim_{(A)} \GradAlg(H)= 46$. To do
  so we use Proposition~\ref{mainlev2}. Indeed from the Hilbert
  functions we see that $s(I_A)=s(I_{A_2})=4$ and $s(I_{A_1})=2$.
  Hence $(I_A \otimes I_{A_2})_7=0$.  Moreover since $I_{A}$ has one
  generator of degree $4$ and $8$ generators of degree $5$ and $A_1$
  is a complete intersection of type $(2,2,6)$, it follows that all
  relations of $I_A \cdot I_{A_1}$ must be of degree greater or equal
  to $8$. We get that $(I_A \otimes I_{A_1})_7 \simeq (I_A \cdot
  I_{A_1})_7$ is an isomorphism of vector spaces of
  dimension $2 \cdot(3+8)=22$.  Hence Proposition~\ref{mainlev2}
  applies and we get the unobstructedness of $A$ and
  $$
  \ \dim_{(A)} \GradAlg(H)= \ 2 \cdot \dim (I_A)_7 - dim (I_A \cdot I_{A_1})_7
  = 46 \ .$$
 
  To get the other component, let now $F_1$ be a sum of length $9$ of general
  linear forms to the $7$-th power (i.e. $H_{A_1}= (1,3,6,9,9,6,3,1)$), let
  $F_2$ be, say $F_2=x^6y+xy^6+z^7$ and let $A=R/ann(F_1,F_2)$. Then the the
  Hilbert function of $A$ is $H$ by Remark~\ref{remcomp}. We claim that $A$ is
  {\it licci}! Indeed it is easily checked by Macaulay 2 that $A$ above admits
  the following CI-linkages to a CI of type $(1,1,3)$. We start with a general
  CI of type $(4,5,7)$ whose generators are contained in $I_A$ and follow up
  by general CI-linkages of type $(4,5,6)$, $(4,4,6)$, $(4,4,5)$, $(3,3,5)$,
  $(3,3,4)$, $(2,2,4)$ and $(2,2,3)$, in this order. Then $A$ is unobstructed
  and $\ \dim_{(A)} \GradAlg(H)=47$ by Remark~\ref{summarize2}(iii) or
  Theorem~\ref{mainlink} and we are done (of course, once using Macaulay 2 it
  is easier to see that the tangent space is 47-dimensional by computing $_0\!
  \ext^1(I_A,I_A)$. The unobstructedness of $A$ is, however, not at all easy
  to see by Macaulay 2 computations because $_0\! \Ext_R^1(I_A,A) \simeq \ 
  _0\!  \Ext^2(I_A,I_A)$ is 1-dimensional and so is $ _0\! \HH^2(R,A,A)$ by
  Proposition~\ref{mainlev2} and Remark~\ref{remmainlev}. Hence we really need
  to use that the unobstructedness property is preserved under CI-linkages,
  which is true by Theorem~\ref{mainlink}).
\end{example}

\begin{remark} \label{remlev2comp} (a) We have tried to look for other
  examples of several ``level type 2 components'' of {\rm smaller} socle
  degree, but have not yet fully succeeded. A promising candidate is
  $H=(1,3,6,9,9,6,2)$ where we get a level type 2 algebra $A$ by
  starting with a CI of type $(2,2,3)$ and linking in one step via a
  CI of type $(4,4,3)$.  By Remark~\ref{summarize2}(iii) we have
  $_0\! \hom_R(I_A,A)=33$.
  Moreover we have an $A'=R/ann(G_1,G_2)$ with $_0\! \hom_R(I_{A'},A')=35$
  (checked by Macaulay 2) by taking $G_1$ (resp. $G_2$) to be a sum of length
  $3$ (resp. $6$) of general linear forms to the $6$-th power. It follows that
  $A$ belongs to a $33$-dimensional generically smooth component while, due to
  the size of the tangent spaces, there are only two possibilities for $A'$.
  It is either obstructed, or it is unobstructed in which case it belongs to an
  irreducible component different from the ``licci'' component. We
  have not yet been able to decide which of the possibilities that occur.
 
  (b) One may construct other examples of several ``level type 2 components''
  of {\rm larger} socle type by taking the two components of
  Example~\ref{lev2comp} and performing a biliaison, starting with general
  CI's of type $(5,5,b)$ containing the general elements of the components and
  follow up by general CI-linkages of type $(b,b,b)$, $b \geq 7$.  Using
  Theorem~\ref{mainlink}, we get two irreducible components of $ \GradAlg(H')$
  whose general elements are level type 2 quotients of socle degree $3b-8$
  ($H'$ may be computed from $H=(1,3,6,10,14,10,6,2)$).
\end{remark}

\begin{remark}
So if we want to compare the parameter space of type 2 codimension 3
level algebras to the corresponding space of Gorenstein algebras, we
see many differences. In the level type 2 case,

(i) \ \ the parameter space may be reducible (Example~\ref{lev2comp} and
Remark~\ref{remlev2comp}(b)),

(ii) \ $ _0\! \HH^2(R,A,A)$ may be non-vanishing
(e.g. Example~\ref{lev2comp}, there are many more). \\
In the Gorenstein case (i) and (ii) are false. We have, however, not yet been
able to find two irreducible ``type 2 codimension 3 level components'' with a
type 2 level algebra in the intersection, nor have we been able to find an
obstructed type 2 codimension 3 level algebra.
\end{remark}

{\small

\bibliographystyle{amsalpha}

} Oslo University College, Faculty of Engineering, Pb.\! 4 St.\! Olavs plass,
N-0130 Oslo, Norway.
 
 E-mail address: JanOddvar.Kleppe@iu.hio.no \hspace*{35.5mm}
 Date: 19. January 2006.

\end{document}